\documentclass[12pt]{amsart}
\usepackage[foot]{amsaddr}
\usepackage{amsmath}
\usepackage{amssymb}
\usepackage{amsthm}
\usepackage{array}
\usepackage[english]{babel}
\usepackage{bigints}
\usepackage{booktabs}
\usepackage{caption}
\usepackage{graphicx, hyperref}
\usepackage{indentfirst}
\usepackage{mathptmx}
\usepackage{mathtools}
\usepackage{multirow}
\usepackage{xcolor}
\usepackage[numbers]{natbib}
\usepackage[paperwidth=8.5in, paperheight=11in, margin=0.5in,headheight=64pt,marginparwidth=0pt,textheight=650pt]{geometry}

\title{Cointegration and unit root tests: a fully Bayesian approach}
\author{Marcio Alves Diniz \\ Federal University of  S\~{a}o Carlos \\
\\
Carlos Alberto de  Bragan\c{c}a Pereira \\ Federal University of Mato Grosso do Sul \\
\\
Julio Michael Stern \\ University of S\~ao Paulo}
\date{}

\begin{document}

\maketitle

\begin{abstract}

To perform statistical inference for time series, one should be able to assess if they present deterministic or stochastic trends.
For univariate analysis, one way to detect stochastic trends is to test if the series has unit roots, 
and for multivariate studies it is often relevant to search  for stationary linear relationships between the series, or if they cointegrate.
The main goal of this article is to briefly review the shortcomings of unit root and cointegration tests proposed by the Bayesian approach of statistical inference and to show how they can be overcome by the Full Bayesian Significance Test (FBST), a procedure designed to test sharp or precise hypothesis.
We will compare its performance with the most used frequentist alternatives, namely, the Augmented Dickey--Fuller for unit roots and the maximum eigenvalue test for cointegration.


\end{abstract}

\section{Introduction}

Several times series present deterministic or stochastic trends, which imply that the effects of these trends on the level of the series are permanent.
Consequently, the mean and variance of the series will not be constant and will not revert to a long-term value.
This feature reflects the fact that the stochastic processes generating these series are not (weakly) stationary, imposing problems to perform inductive inference using the most traditional estimators or predictors.
This is so because the usual properties of these procedures will not be valid under such conditions.

Therefore, when modeling non-stationary time series, one should be able to properly detrend the used series, either by directly modeling the trend by deterministic functions, or by transforming the series to remove stochastic trends.
To determine which strategy is the suitable solution, several statistical tests were developed since the 1970s by the frequentist school of statistical inference.

The Augmented Dickey--Fuller (ADF) test is one of the most popular tests used to assess if a time series has a stochastic trend or, for series described by auto-regressive models, if they have a unit root.
When one is searching for long term relationships between multiple series under analysis, it is crucial to know if there are stationary linear combinations of these series, i.e., if the series are cointegrated.
Cointegration tests were developed, also by the frequentist school, in the late 1980s \citep{engle1987} and early 1990s~\cite{johansen1990}.
Only in the late 1980s did the Bayesian approach to test the presence of unit roots start to be developed.

Both unit root and cointegration tests may be considered tests on precise or sharp hypotheses, i.e., those in which the dimension of the parameter space under the tested hypothesis is smaller than the dimension of the unrestricted parameter space.
Testing sharp hypotheses poses major difficulties for either the frequentist or Bayesian paradigms, such as the need to eliminate nuisance parameters.
 
The main goal of this article is to briefly review the shortcomings of the tests proposed by the Bayesian school and how they can be overcome by the Full Bayesian Significance Test (FBST).
More specifically, we will compare its performance with the most used frequentist alternatives, the ADF for unit roots, and the maximum eigenvalue test for cointegration. 
Since this is a review article, it is important to remark that the results presented here were published elsewhere by the same authors, see \cite{diniz2011,diniz2012a}.

To accomplish this objective, we will define the FBST in the next section, also showing how it can be implemented in a general context.
The following section discusses the problems of testing the existence of unit roots in univariate time series and how the Bayesian tests approach the problem.
Section 4 then shows how the FBST is applied to test if a time series has unit roots and illustrates this with applications on a real data set.
In the sequel, we discuss the Bayesian alternatives to cointegration tests and then apply the FBST to test for cointegration using real data sets.
We conclude with some remarks and possible extensions for future work.

\section{FBST} 
\label{sec:FBST}

The Full Bayesian Significance Test was proposed in \cite{pereira1999} mainly to deal with sharp hypotheses.
The procedure has several properties, see \cite{pereira2008, stern2020}, most interestingly the fact that it is only based on posterior densities, thus avoiding the necessity of complications such as the elimination of nuisance parameters or the adoption of priors with positive probabilities attached to sets of zero Lebesgue measure.

We shall consider general statistical models in which the parameter space is denoted by $\Theta \subseteq \mathbb{R}^m$, $m\in\mathbb{N}$. 
A sharp hypothesis $H$ assumes that $\theta$, the parameter vector of the chosen statistical model, belongs to a sub-manifold $\Theta_H$ of smaller dimensions than $\Theta$.
This implies, for continuous parameter spaces, that the subset $\Theta_H$ has null Lebesgue measure whenever $H$ is sharp.
The sample space, the set of all possible values of the observable random variables (or vectors), is here denoted by $\mathcal{X}$.

Following the Bayesian paradigm, let $h(\cdot)$ be a probability prior density over $\Theta$, ${\bf x}\in\mathcal{X}$, the observed sample (scalar or vector), and $L(\cdot\mid{\bf x})$ the likelihood derived from data ${\bf x}$. 
To evaluate the Bayesian evidence based on the FBST, the sole relevant entity is the posterior probability density for $\theta$ given ${\bf x}$, 

\begin{displaymath}
g(\theta\mid{\bf x}) \propto h(\theta) \cdot L(\theta\mid{\bf x}).
\end{displaymath}

It is important to highlight that the procedure may be used when the parameter space is discrete.
However, when the posterior probability distribution over $\Theta$ is absolutely continuous, the FBST appears as a more suitable alternative to significance hypothesis testing. 
For notational simplicity, we will denote $\Theta_H$ by $H$ in the sequel.

Let $r(\theta)$ be a reference density on $\Theta$ such that the function $s(\theta)=g(\theta\mid{\bf x})/r(\theta)$ is a {\it relative surprise}, (see \cite{good1983}, p.145--146) function.
The reference density is important because it guarantees that the FBST is invariant to reparametrizations, even when $r(\theta)$ is improper, see \cite{pereira2008, stern2013}.
Thus, when considering $r(\theta)$ proportional to a constant, the surprise function will be, in practical terms, equivalent to the posterior distribution.  
For the applications considered in this article, we will use the improper uniform density as reference density on $\Theta$. 
The authors of \cite{madruga2001} remark that it is possible to generalize the procedure using other reference densities such as neutral, invariant, maximum-entropy or non-informative priors, if they are available and desirable.

\

\noindent{\textbf{Definition 1}} [{\bf Tangent set}]. \emph{Considering a sharp hypothesis $H: \theta \in \Theta_H$, the tangential set of the hypothesis given the sample is given by}

\begin{equation}
\mathbb{T}_{{\bf x}}=\{\theta \in \Theta: s(\theta) > s^*\}.
\end{equation}

\noindent
\emph{where} $s^*=\sup_{\theta \in H} s(\theta)$.\\

Notice that the tangent set $\mathbb{T}_{{\bf x}}$ is the highest relative surprise set, that is, the set of points of the parameter space with higher relative surprise than any point in $H$, being {\it tangential} to $H$ in this sense.
This approach takes into consideration the statistical model in which the hypothesis is defined, using several components of the model to define an evidential measure favoring the hypothesis. 
    
\

\noindent{\textbf{Definition 2}} [{\bf Evidence}]. \emph{The Bayesian evidence value against $H$, $\overline{ev}$, is defined as}

\begin{equation}\label{eq:evc}
\overline{ev}=P\left(\theta \in \mathbb{T}_{{\bf x}} \mid {\bf x}\right)=\int\limits_{\mathbb{T}_{{\bf x}}} dG_{\bf x}(\theta),    
\end{equation}

\noindent
\emph{where $G_{\bf x}(\theta)$ denotes the posterior distribution function of $\theta$ and the above integral is of the Riemann--Stieltjes~type.}

\

Definition 2 sets $\overline{ev}$ as the posterior probability of the tangent set that is interpreted as an evidence value against $H$. 
Hence, the evidence value supporting $H$ is the complement of $\overline{ev}$, namely, $ev=1-\overline{ev}$.
Notwithstanding, $ev$ is not evidence against $A: \theta\notin \Theta_H$, the alternative hypothesis (which is not sharp anyway). Equivalently, $\overline{ev}$ is not evidence in favor of $A$, although it is against $H$.

\

\noindent{\textbf{Definition 3}} [{\bf Test}]. \emph{The FBST is the procedure that rejects $H$ whenever $ev=1-\overline{ev}$ is smaller than a critical level, $ev_c$.}    

\

Thus, we are left with the problem of deciding the critical level $ev_c$ for each particular application.
We briefly discuss this and other practical issues in the following subsection.

\subsection{Practical Implementation: Critical Values and Numerical Computation}
\label{subsec:practical}

Since $ev$ (also called e-value) is a statistic, it has a sampling distribution derived from the adopted statistical model and in principle this distribution could be used to find a threshold value.
If the likelihood and the posterior distribution satisfy certain regularity conditions. See \cite{schervish1995}, p. 436.
 \cite{diniz2012b} proved that, asymptotically, there is a relationship between $ev$ and the $p$-values obtained from the frequentist likelihood ratio procedure used to test the same hypotheses.
This fact provides a way to find, at least asymptotically, a critical value to $ev$ to reject the hypothesis being tested.

In a recent review \cite{stern2020}, the authors discuss different ways to provide a threshold for $ev$.
Among these alternatives, we highlight the standardized e-value, which follows, asymptotically, the uniform distribution on $(0,1)$. See also \cite{borges2007} for more on the standardized version of $ev$. 

One could also try to define the FBST as a Bayes test derived from a particular loss function and the respective minimization of the posterior expected loss.
Following this strategy, \cite{madruga2001} showed that there are loss functions which result in $ev$ as a Bayes estimator of $\phi=\mathbb{I}_H(\theta)$, where $\mathbb{I}_A(x)$ denotes the indicator function, being equal to one if $x\in A$ and zero otherwise, $x\notin A$. 
Hence, the FBST is in fact a Bayes procedure in the formal sense as defined by Wald in \cite{wald1950}.

\begin{table}[!h]
\centering
\caption{Pseudocode to implement the FBST}
\begin{tabular}{l}
\hline
\textbf{General algorithm}: compute $ev$ supporting hypothesis $H: \theta\in\Theta_H$\\
\hline
1. Specify the statistical model (likelihood) and prior distribution on $\Theta$.\\
2. Specify the reference density, $r(\theta)$, and derive the relative surprise function, $s(\theta)$.\\
3. Find $s^*$, the maximum value of $s(\theta)$ under the constraint $\theta\in H$.\\
4. Integrate the posterior distribution on the tangent set---Equation (2)---to find $\overline{ev}$.\\
5. Find $ev=1-\overline{ev}$.\\
\hline
\end{tabular}
\label{tab:pseudocode}
\end{table}

To compute the evidence value supporting $H$ defined in the last section, we need to follow the steps showed in Table \ref{tab:pseudocode}. Appendix \ref{app:computer} provides detailed information about the computational resources and codes used to implement the FBST in the examples presented in this work.
After defining the statistical model and prior, it is simple to find the surprise function, $s(\theta)$.
In step 3, one should find the point of the parameter space  in $H$ that maximizes $s(\theta)$, that is, to solve a problem of  constrained numerical maximization.
In several applications, this step does not present a closed form solution, requiring the use of numerical optimizers.

Step 4 involves the integration of the posterior distribution on a subset of $\Theta$, the tangent set $\mathbb{T}_{{\bf x}}$ that can be highly complex.
Once more, since in many cases it is fairly difficult to find an explicit expression for $\mathbb{T}_{{\bf x}}$, one may use various numerical techniques to compute the integral. 
If it is possible to generate random samples from the posterior distribution, Monte Carlo integration provides an estimate of $ev$, as we will show in this work.
Another alternative is to use approximation techniques, such as those proposed in \cite{kadane1986}, based on a Laplace approximation.
We discuss how to implement such approximations for unit root and cointegration tests in \cite{diniz2011,diniz2012a}.

\section{Bayesian Unit Root Tests} \label{tests}

Before presenting the Bayesian procedures used to test the presence of unit roots, let us fix notation.
We will denote by $y_t$ the $t$-th value of a univariate time series observed in $t=1,\ldots,T+p$ dates, where $T$ and $p$ are positive integers.
The usual approach is to assume that the series under analysis is described by an auto-regressive  process with $p$ lags, $AR(p)$, 
meaning that the data generating process is fully described by a stochastic difference equation of order $p$, possibly with an intercept or drift and a deterministic linear trend, i.e.,

\begin{equation}\label{eq:ar_p}
y_t=\mu+\delta \cdot t+\phi_1 y_{t-1}+\ldots+\phi_p y_{t-p}+\varepsilon_t
\end{equation}

\noindent
with $\varepsilon_t$ i.i.d. $N(0,\sigma^2)$ for $t=1,\ldots,T+p$.
Using the lag or backshift operator $B$, we denote $y_{t-k}$ by $B^ky_t$, allowing us to rewrite \eqref{eq:ar_p} as

\begin{equation}
(1-\phi_1B-\ldots \phi_pB^p)y_t=\mu+\delta \cdot t+\varepsilon_t
\end{equation}

\noindent
where $\phi(B)=(1-\phi_1B-\ldots \phi_pB^p)$ is the autoregressive polynomial.
The difference Equation \eqref{eq:ar_p} will be stable, implying that the process generating $\{y_t\}_{t=1}^{T+p}$ is (weakly) stationary, whenever the roots of the characteristic polynomial $\phi(z)$, $z\in\mathbb{C}$, lie outside the unit circle, since there may be complex roots. The set of polynomial operators, such as lag polynomials like $\phi(B)$, induces an algebra that is isomorphic to the algebra of polynomials in real or complex variables, see \cite{dhrymes2000}.

If some of the roots lie exactly on the unit circle, it is said that the process has unit roots.
In order to test such a hypothesis statistically, \eqref{eq:ar_p} is rewritten as

\begin{equation}\label{eq:ecm}
\Delta y_t=\mu+\delta \cdot t+\Gamma_0 ~ y_{t-1}+\Gamma_1 \Delta y_{t-1}+\ldots+\Gamma_{p-1} \Delta y_{t-p+1}+\varepsilon_t
\end{equation}

\noindent
where $\Delta y_t=y_t-y_{t-1}$, $\Gamma_0=\phi_1+\ldots+\phi_p-1$ and $\Gamma_i=-\sum_{j=i+1}^{p}\phi_j$, for $i=1,\ldots,p-1$. 
If the generating process has only one unit root, one root of the complex polynomial $\phi(z)$, 

\[
1-\phi_1z-\phi_2z^2\ldots \phi_pz^p,
\]

\noindent
is equal to one, meaning that

\[
1-\phi_1-\phi_2-\ldots-\phi_p=0
\]

\noindent
i.e., $\phi(1)=0$, and all the other roots are on or outside the unit circle.
In this case, $\Gamma_0=0$, the hypothesis that will be tested when modeling \eqref{eq:ecm}.
Even though tests based on these assumptions verify if the process has a single unit root, there are generalizations based on the same principles that test the existence of multiple unit roots, see \cite{dickey1987}.

The search for Bayesian unit root tests began in the late 1980s. 
As far as we know, \cite{sims1988} and \cite{sims1991} were the first works to propose a Bayesian approach for unit root tests.
The frequentist critics of these articles received a proper answer in \cite{phillips1991a,phillips1991b}, generating a fruitful debate that produced a long list of papers in the literature of Bayesian time series. 
A good summary of the debate and the Bayesian papers that resulted from it is presented in \cite{bauwens1999}.
We will present here only the most relevant strategies proposed by the Bayesian school to test for unit roots.

Let $\theta=(\rho,\psi)$ be the parameters vector, in which $\rho=\sum_{i=1}^p \phi_i$ and $\psi=(\mu,\delta,\Gamma_1,\ldots,\Gamma_{p-1})$.
Assuming $\sigma^2$ fixed, the prior density for $\theta$ can be factorized as

\begin{displaymath}
h(\theta)=h_0(\rho)\cdot h_1(\psi\mid\rho).
\end{displaymath}

The marginal likelihood for $\rho$, denoted by $L_m$, is:

\begin{displaymath}
L_m(\rho\mid{\bf y})\propto \int\limits_{\Psi} L(\theta\mid{\bf y})\cdot h_1(\psi\mid\rho) ~ d\psi.
\end{displaymath}

\noindent
where ${\bf y}=\{y_t\}_{t=1}^{T+p}$ is the observations vector, $L(\theta|{\bf y})$ the full likelihood, and $\Psi$ the support of the random vector $\psi$. 
This marginal likelihood, associated with a prior for $\rho$, is the main ingredient used by standard Bayesian procedures to test the existence of unit roots. Even though the procedure varies among authors according to some specific aspects, mentioned below, basically all of them use Bayes factors and posterior probabilities.

One important issue is the specification of the null hypothesis: some authors, starting from \cite{schotman1991}, consider $H_0:\rho=1$ against $H_1:\rho<1$. 
Starting from \cite{dickey1979}, this is the way the frequentist school addresses the problem, but following this approach no explosive value for $\rho$ is considered. 
The decision theoretic Bayesian approach solved the problem using the posterior probabilities ratio or Bayes factor:

\begin{displaymath}
B_{01}=\frac{L_m(\rho=1\mid{\bf y})}{\int\limits_0^1L_m(\rho\mid{\bf y})\cdot h_0(\rho) ~ d\rho}.
\end{displaymath}

Advocates of this solution argue that one of the advantages of this approach is that the null and the alternative hypotheses are given equal weight.
However, the expression above is not defined if $h_0(\rho)$ is not a proper density since the denominator of the Bayes factor is equal to the predictive density, defined just if $h_0(\rho)$ is a proper density. 
There are also problems if $L_m(\rho=1|{\bf y})$ is zero or infinite.

The problem is approached by \cite{phillips1991a,lubrano1995} by testing  $H_0:\rho\geq1$ against $H_1:\rho<1$, considering explicitly explosive values for $\rho$. 
The main advantage of this strategy is the possibility to compute posterior probabilities like

\begin{displaymath}
P(\rho>1\mid{\bf y})=\int_{1}^{\infty}g_m(\rho\mid{\bf y}) ~ d\rho
\end{displaymath}

\noindent
defined even for improper priors on $\rho$, where $g_m$ is the marginal posterior for $\rho$. 

In \cite{dejong1991}, the authors do not choose $\rho$ as the parameter of interest, examining instead the largest absolute value of the roots of the characteristic polynomial and then verifying if it is smaller or larger than one. 
Usually, this value is slightly smaller than $\rho$, but the authors argue that this small difference may be important. 
When this approach is used, unit roots are found less frequently. For an AR(3) model with a constant and deterministic trend, \cite{dejong1991} derives the posterior density for the dominant root for the 14 series used in \cite{nelson1982} and concluded the following: for eleven of the series, the dominant root was smaller than one, that is to say, the series were trend-stationary. 
These results were based on flat priors for the autoregressive parameters and the deterministic trend coefficient.

Another controversy is about the prior over $\rho$: \cite{phillips1991a} argues that the difference between the results given by the frequentist and Bayesian inferences is due to the fact that the flat prior proposed in \cite{sims1988} overweights the stationary region of $\rho$. Hence, he derived a Jeffreys prior for the AR(1) model: this prior quickly diverges as $\rho$ increases and becomes larger than one. 
The obtained posterior led to the same results of \cite{nelson1982}, which will be discussed in detail in the following  section. 
The critics of the approach adopted by Phillips in \cite{phillips1991a} judged the Jeffreys prior as unrealistic, from a subjective point of view.
See the comments on Phillips's paper on the \textit{Journal of Applied Econometrics}, volume 6, number 4, 1991. The subsequent papers of the same number support the Bayesian approach.
This is a nonsensical objection if one considers that the Jeffreys prior is crucial to ensure an invariant inferential procedure, and invariance is a highly desirable property, for either objective or subjective reasons.
See \cite{stern2011} for more on invariance in physics and statistical models. 

A final controversial point concerns the modeling of initial observations.
If the likelihood explicitly models the initial observed values (it is an {\it exact} likelihood), the process is implicitly considered  stationary. 
In fact, when it is known that the process is stationary, and it is believed that the data generating process is working for a long period, it is reasonable to assume that the  parameters of the model determine the marginal distribution of the initial observations. 
In the simplest AR(1) model, this would imply that $y_1 \sim N(0,\sigma^2/(1-\rho^2))$. 
In this scenario, to perform the inference conditional on the first observation would discard relevant information. 
On the other hand, there is no marginal distribution defined for $y_1$ if the generating process is not stationary. 
Then, it is valid to use a likelihood conditional on initial observations. 
For the models presented here, we always work with the conditional likelihood.
As argued in \cite{sims1988}, inferences for stationary models are little affected by using conditional likelihoods, especially for large samples. 
He compares these inferences with the ones based on exact likelihoods under explicit modeling for initial observations.

\section{Implementing the FBST for Unit Root Testing}
\label{sec:FBST_unit_root}

We will now describe how to use the FBST to test for the presence of unit roots referring to the general model \eqref{eq:ecm}.
It is also possible to include $q\in\mathbb{N}$ moving average terms in \eqref{eq:ar_p} to model the process, a case that will not be covered in this article but that, in principle, shall not imply major problems for the FBST.

\begin{equation}
\Delta y_t=\mu+\delta \cdot t+\Gamma_0 ~ y_{t-1}+\Gamma_1 \Delta y_{t-1}+\ldots+\Gamma_{p-1} \Delta y_{t-p+1}+\varepsilon_t, \tag{\ref{eq:ecm}}
\end{equation}

\noindent
where $\varepsilon_t\stackrel{i.i.d.}\sim N(0,\sigma^2)$ for $t=1,\ldots,T+p$, recalling also that the hypothesis being tested is $\Gamma_0=0$.
We slightly change the notation of the last section now using $\psi$ to denote the vector $(\mu, \delta, \Gamma_0, \ldots, \Gamma_{p-1})$ and setting $\theta=(\psi, \sigma)$.

Recalling the steps to implement the FBST displayed in Table \ref{tab:pseudocode}, we have just specified the statistical model.
The likelihood, conditional on the first $p$ observations, derived from the Gaussian model is  

\begin{equation}
L(\theta\mid{\bf y}) = (2 \pi)^{-T/2}\sigma^{-T}\textrm{exp}\left\{- \frac{1}{2 \sigma^2}\cdot\sum_{t=p+1}^{T+p}\varepsilon_t^2 \right\},
\end{equation}

\noindent{in which $\varepsilon_t=\Delta y_t-\mu-\delta \cdot t-\Gamma_0y_{t-1}-\Gamma_1\Delta y_{t-1}-\ldots-\Gamma_{p-1}\Delta y_{t-p+1}$.}
To complete step 1 of Table \ref{tab:pseudocode}, we need a prior distribution for $\theta$.
For all the series modeled in this article, we will use the following non informative prior:

\begin{equation}
h(\theta)=h(\psi,\sigma) \propto 1/\sigma.
\end{equation}

We are aware of the problems caused by improper priors applied to this problem when one uses alternative approaches, like those mentioned by \cite{bauwens1999}.
However, one of our goals is to show how the FBST can be implemented even for a potentially problematic prior like this one. 
To write the posterior, we use the following notation:

\[
\Delta Y=\left[
\begin{array}{c}
\Delta y_{p+1} \\
\Delta y_{p+2}\\
\vdots\\
\Delta y_{T+p}\\
\end{array}
\right], ~ ~ 
X=\left[ \begin{array}{cccccc}
  1   &p+1   &y_p   & \Delta y_p  &\ldots   &\Delta y_2 \\
  1   &p+2   &y_{p+1} & \Delta y_{p+1}     &\ldots   &\Delta y_3\\
  \vdots&\vdots&\vdots&\vdots&\vdots&\vdots\\
  1   &T+p   &y_{T+p-1} &\Delta y_{T+p-1}   &\ldots   &\Delta y_{T+1}\\
  \end{array} \right], ~ ~
\psi=\left[
\begin{array}{c}
\mu \\
\delta\\
\Gamma_0\\
\vdots\\
\Gamma_{p-1}\\
\end{array}
\right],
\]

\noindent{being $\Delta Y$ of dimension $T \times 1$, $X$ of dimension $T \times (p+2)$ and $\psi$, $(p+2) \times 1$.}
Thanks to this notation, we can write, using primes to denote transposed matrices:

\[
\sum_{t=p+1}^{T+p}\varepsilon_t^2=(\Delta Y-X\psi)'(\Delta Y-X\psi)=(\Delta Y-\widehat{\Delta Y})'(\Delta Y-\widehat{\Delta Y})+(\psi-\widehat{\psi})'X'X(\psi-\widehat{\psi}),
\]

\noindent
where $\widehat{\psi}=(X'X)^{-1}X'\cdot\Delta Y$ is the ordinary least squares (OLS) estimator of $\psi$ and $\widehat{\Delta Y}=X\widehat{\psi}$ its prediction for $\Delta Y$.
Thus, the full posterior is

\begin{equation}
g(\theta\mid{\bf y}) \propto \sigma^{-(T+1)} \textrm{exp}\left\{-\frac{1}{2\sigma^2}[(\Delta Y-\widehat{\Delta Y})'(\Delta Y-\widehat{\Delta Y})+(\psi-\widehat{\psi})'X'X(\psi-\widehat{\psi})]\right\}, \label{eq:post}
\end{equation}

\noindent
a Normal-Inverse Gamma density.

Step 2 demands a reference density in order to define the relative surprise function.
Since we will use the improper density $r(\theta)\propto 1$, the surprise function will be equivalent to the posterior distribution in our applications.
Given this, to find $s^*$ (Step 3), we need to find the maximum value of the posterior under the hypothesis being tested, in our case, $\Gamma_0=0$.

This maximization step is fairly simple to implement given the modeling choices made here: Gaussian likelihood, non informative prior and reference density proportional to a constant.
The restricted (assuming $H$) posterior distribution is

\begin{equation}
g_r(\theta_r\mid{\bf y}) \propto \sigma^{-(T+1)} \textrm{exp}\left\{-\frac{1}{2\sigma^2}[(\Delta Y-\widehat{\Delta Y}_r)'(\Delta Y-\widehat{\Delta Y}_r)+(\psi_r-\widehat{\psi}_r)'X_r'X_r(\psi_r-\widehat{\psi}_r)]\right\}, \label{eq:post_r}
\end{equation}

\noindent
in which $\theta_r=(\psi_r,\sigma)$, $\psi_r$ being vector $\psi $ without $\Gamma_0$,

\[
X_r=\left[ \begin{array}{ccccc}
  1   &p+1   & \Delta y_p  &\ldots   &\Delta y_2 \\
  1   &p+2  & \Delta y_{p+1}     &\ldots   &\Delta y_3\\
  \vdots&\vdots&\vdots&\vdots&\vdots\\
  1   &T+p  &\Delta y_{T+p-1}   &\ldots   &\Delta y_{T+1}\\
  \end{array} \right], ~ ~
\widehat{\psi}_r=(X_r'X_r)^{-1}X_r'\cdot\Delta Y, ~ \text{and} ~ ~
\widehat{\Delta Y}_r=X_r\widehat{\psi}_r,
\]

\noindent
that is, $X_r$ is simply matrix $X$ above without its third column, since under $H:\Gamma_0=0$ and the coefficient of the third column of $X$ is $\Gamma_0$,--see Equation \eqref{eq:ecm}--$\widehat{\psi}_r$ is a least squares estimator of $\psi_r$ and $\widehat{\Delta Y}_r$ denotes the predicted values for $\Delta Y$ given by the restricted model.
From \eqref{eq:post_r}, it is easy to show that the maximum a posteriori (MAP) estimator of $\theta_r$ is given by $(\widehat{\psi}_r,\widehat{\sigma}_r)$, with

\[
\widehat{\sigma}_r=\sqrt{\frac{(\Delta Y-\widehat{\Delta Y}_r)'(\Delta Y-\widehat{\Delta Y}_r)}{T+1}}.
\]

Plugging the values of $\widehat{\psi}_r$ and $\widehat{\sigma}_r$ into \eqref{eq:post_r}, we find $s^*$, as requested in Step 3.
Step 4 will also be easy to implement thanks to structure of the models assumed in this section.
Since the full posterior, \eqref{eq:post}, is a Normal-Inverse Gamma density, a simple Gibbs sampler allows us to obtain a random sample from such distribution, suggesting a Monte Carlo approach to compute $\overline{ev}$.
From \eqref{eq:post}, the conditional posteriors of $\psi$ and $\sigma$ are, respectively,
 
\begin{align}
g_{\psi}(\psi\mid\sigma,{\bf y})&\propto
  N(\widehat{\psi},\sigma^2(X'X)^{-1})\label{eq:cond1}\\
g_{\sigma^2}(\sigma^2\mid\psi,{\bf y})& \propto
 IG\left({T+1\over 2},H\right)\label{eq:cond2}
\end{align}

\noindent
in which $H=0.5[(\Delta Y-\widehat{\Delta Y})'(\Delta Y-\widehat{\Delta Y})+(\psi-\widehat{\psi})'X'X(\psi-\widehat{\psi})]$, $IG$ denotes the Inverse-Gamma distribution and $\widehat{\psi}$ is the OLS estimator of $\psi$, as above. 
Appendix \ref{app:dist} brings the parametrization and the probability density function of the Inverse-Gamma distribution.
With a sizable random sample from the full posterior, we estimate $\overline{ev}$ as the proportion of sampled vectors that generate a value for the posterior greater than $s^*$, found in Step 3.
Hence, in Step 5, we only compute one minus the estimate of $\overline{ev}$ found in Step 4.
The whole procedure is summarized in Table \ref{tab:unit_root_tests}.
For the implementations in this article we sampled 51,000 vectors from \eqref{eq:post} and discarded the first 1,000 as a burn-in sample.

\begin{table}[!h]
\centering
\caption{Pseudocode to implement the FBST to unit root tests}
\begin{tabular}{l}
\hline
\textbf{General algorithm}: compute $ev$ supporting hypothesis $H: \Gamma_0=0$ in model \eqref{eq:ecm}\\
\hline
1. Statistical model: Gaussian; prior: $h(\theta)\propto1/\sigma$.\\
2. Reference density: $r(\theta)\propto 1$; relative surprise function: $g(\theta\mid{\bf y})$.\\
3. Find $s^*$: \eqref{eq:post_r} evaluated at $\widehat{\psi}_r$ and $\widehat{\sigma}_r$.\\
4. Gibbs sampler (from Equations \eqref{eq:cond1} and \eqref{eq:cond2}) to obtain $N$ random samples of parameter vectors from \eqref{eq:post}.\\
~ ~ Evaluate the posterior at the sampled vectors and estimate $\overline{ev}$ as the proportion of $N$ in which the \\
~ ~ evaluated values are larger than $s^*$.\\
5. Find $ev=1-\overline{ev}$.\\
\hline
\end{tabular}

\label{tab:unit_root_tests}
\end{table}

\subsection{Results}

We implemented the FBST as described above to test the presence of unit roots in 14 U.S. macroeconomic time series, all with annual frequency, first mentioned in \cite{nelson1982}.
We used the extended series, analyzed in \cite{schotman1991}.
Appendix \ref{app:computer} brings more information on the data set and the computational resources and codes used to obtain the results displayed in Table \ref{tab:ur_results} below.

Table \ref{tab:ur_results} reports the names of the tested series, the number of available observations or sample size, the adopted value for $p$---as denoted in Equation \eqref{eq:post}---, if a linear (deterministic) trend was included in the model or not, the ADF test statistic and its respective $p$-value.
We have used the computer package described in \cite{mackinnon1994} to find the ADF $p$-values, available in the R library \texttt{urca}.
The last two columns bring the posterior probability of non-stationarity, $\Gamma_0\geq0$, and the FBST e-values for the specified models.
In order to obtain comparable results, we have adopted the models chosen by \cite{bauwens1999} for all the series.
All the models considered the intercept or constant term, $\mu$ in \eqref{eq:post}.

The results show that the non-stationary posterior probabilities are quite distant from the ADF $p$-values. 
These results were highlighted in \cite{sims1988,sims1991}. 
Considering the simplest AR(1) model, they argued that, once frequentist inference is based on the distribution of $\widehat{\rho}|\rho=1$, the non-stationary posterior probabilities  provide counterintuitive conclusions since the referred distribution is skewed. 
Their main argument is that Bayesian inference uses a distribution (marginal posterior of $\rho$) that is not skewed. 

As mentioned before, ref. \cite{phillips1991a} claims that the difference in results between frequentist and Bayesian approaches is due to the flat prior that puts much weight on the stationary region. 
He proposed the use of Jeffreys priors, which restored the conclusions drawn by the frequentist test. 
Phillips argued that the flat prior was, actually, informative when used in time series models like those for unit root tests. 
Using simulations, he shows that
{\it `` [the use of a] flat prior has a tendency to bias the posterior  towards stationarity. ...  even when [the estimate] is close to unity, there may still be a non negligible downward bias in the [flat] posterior
 probabilities''}.
Notwithstanding, the e-values reported in the last column are quite close to the ADF $p$-values even using the flat prior criticized by Phillips.

\begin{table}[!h]
\centering
\caption{Unit root tests for the extended Nelson and Plosser data set}\label{tab:ur_results}
\begin{tabular}{lccccccc}
\toprule
\textbf{Series }				&\textbf{Sample Size}  	&\boldmath{$p$}	&\textbf{Trend	}	&\textbf{ADF}	&\boldmath{$p$}\textbf{-Value}	&\boldmath{$P(\Gamma_0\geq0|{\bf y})$}& \textbf{e-Value} \\ 
\midrule
Real GNP				&80		&2	&yes	 		&$-3.52$	&0.044	&0.0005&	0.040	\\
Nominal GNP				&80		&2	&yes			&$-2.06$	&0.559	&0.0238&	0.523	\\
Real GNP per capita		&80		&2	&yes			&$-3.59$	&0.037	&0.0004&	0.034	\\
Industrial prod.		&129	&2	&yes			&$-3.62$	&0.032	&0.0003&	0.028	\\
Employment				&99 	&2	&yes			&$-3.47$	&0.048	&0.0004&	0.043	\\
Unemployment rate		&99 	&4	&no				&$-4.04$	&0.019	&0.0001&	0.020	\\
GNP deflator			&100	&2	&yes			&$-1.62$	&0.778	&0.0584&	0.762	\\
Consumer prices			&129	&4	&yes			&$-1.22$	&0.902	&0.1154&	0.983	\\
Nominal wages			&89 	&2	&yes			&$-2.40$	&0.377	&0.0106&	0.341	\\
Real	wages			&89 	&2	&yes			&$-1.71$	&0.739	&0.0475&	0.715	\\
Money stock			 	&100 	&2	&yes			&$-2.91$	&0.164	&0.0029&	0.147	\\
Velocity				&119	&2	&yes			&$-1.62$	&0.779	&0.0620&	0.777	\\
Bond yield				&89 	&4	&no				&$-1.35$	&0.602	&0.0962&	0.936	\\
Stock prices			&118	&2	&yes			&$-2.44$	&0.357	&0.0103&	0.349	\\ \hline
\end{tabular}
\end{table}

\section{Bayesian Cointegration Tests}
\label{sec:bayes_coint}

Before starting our brief review of the most relevant Bayesian cointegration tests, we fix notation and present the definitions to which we will refer in the sequel.

All the tests mentioned here are based on the following multivariate framework.
Let ${\bf Y}_t=[y_{1t} \ldots y_{nt}]'$ be a vector with $n\in\mathbb{N}$ time series, all of them assumed to be integrated of order $d\in\mathbb{N}$, i.e., have $d$ unit roots. The series are said to be cointegrated if there is a nontrivial linear combination of them that has $b\in\mathbb{N}$ unit roots, $b<d$.
We will assume that, as in most applications, $d=1$ and $b=0$, meaning that, if the time series in ${\bf Y}_t$ is cointegrated, there is a linear combination ${\bf a}'{\bf Y}_t$ that is stationary, where ${\bf  a}\in\mathbb{R}^n$ is the cointegrating vector.
Since the linear combination ${\bf a}' {\bf Y}_t$ is often motivated by problems found in economics, it is called a long-run equilibrium relationship. 
The explanation is that non-stationary time series that are related by a long-run  relationship cannot drift too far from the equilibrium because economic forces will act to restore the relationship.

Notice also that: (i) the cointegrating vector is not uniquely determined since, for any scalar $s$, $(s\cdot{\bf a})$ is a cointegrating vector; and (ii) if ${\bf Y}_t$ has more than two series, it is possible that there is more than one cointegrating vector generating a stationary linear combination.

It is assumed that the data generating process of ${\bf Y}_t$ is described by the following vector autoregression with $p\in\mathbb{N}$ lags, denoted VAR(p), and given by:

\begin{equation}\label{eq:varp}
{\bf Y}_t={\bf c}+\Phi_0 {\bf D}_t+\Phi_1 {\bf Y}_{t-1}+\ldots+\Phi_p {\bf Y}_{t-p}+{\bf E}_t,
\end{equation}

\noindent
in which ${\bf c}$ is a $(n\times 1)$ vector of constants, ${\bf D}_t$ a vector $(n\times 1)$ with some deterministic variable, such as deterministic trends or seasonal dummies, $\Phi_i$ are $(n\times n)$ coefficients matrices and  ${\bf E}_t$ is a $(n\times 1)$ stochastic vector with multivariate normal distribution with null expected value and covariance matrix $\Omega$, denoted $N_n({\bf 0},\Omega)$.
This dynamic model is assumed valid for $t=1,\ldots, T+p$, the available span of observations of ${\bf Y}_t$.
As in the univariate case, one may include moving average terms in \eqref{eq:varp}, i.e., lags for ${\bf E}_t$, but this, in principle, would not cause major problems in the Bayesian framework.
Model \eqref{eq:varp} can be rewritten using the lag or backshift operator as 

\begin{equation}
(I_n-\Phi_1B-\ldots-\Phi_pB^p){\bf Y}_t={\bf c}+\Phi_0 {\bf D}_t+{\bf E}_t,
\end{equation}

\noindent
where $\Phi(B)=I_n-\Phi_1B-\ldots-\Phi_pB^p$ is the (multivariate) autoregressive polynomial and $I_n$ denotes the $n$-dimensional identity matrix.
The associate characteristic polynomial in this context will be the determinant of $\Phi(z)$, $z\in\mathbb{C}$.
If all the roots of the characteristic polynomial lie outside the unit circle, it is possible to show that ${\bf Y}_t$ has a stationary representation---see \cite{johansen1996}---such as Equation \eqref{eq:varp}.
In order to determine if this is the case, model \eqref{eq:varp} is rewritten as an (vectorial) error correction model (VECM):

\begin{equation}
\Delta {\bf Y}_t={\bf c}+\Phi_0 {\bf D}_t+\Gamma_1 \Delta {\bf Y}_{t-1}+\ldots+\Gamma_{p-1}\Delta {\bf Y}_{t-p+1}+\Pi {\bf Y}_{t-1}+{\bf E}_t, \label{eq:vec}
\end{equation}

\noindent
where $\Delta {\bf Y}_t=[\Delta y_{1t}\ldots\Delta y_{nt}]'$, $\Gamma_i=-(\Phi_{i+1}+\ldots\Phi_p)$ for $i=1,2,\ldots,p-1$ and $\Pi=-\Phi(1)=-(I_n-\Phi_1-\ldots-\Phi_p)$.
It is possible to show that, when all the roots of det($\Phi(z)$) are outside the unit circle, matrix $\Pi$ in \eqref{eq:vec} has full rank, i.e., all the $n$ eigenvalues of $\Pi$ are $n$ non null.
If the rank of $\Pi$ is null, this matrix cannot be distinguished from a null matrix, implying that the series in ${\bf Y}_t$ has at least one unit root and a valid representation is a VAR of order $p-1$, i.e., model \eqref{eq:vec} without the term $\Pi{\bf Y}_{t-1}$.
It is possible that the series in ${\bf Y}_t$ has two unit roots each, implying that the correct VECM must be written with $\Delta^2 {\bf Y}_t$ as a dependent variable.

Finally, if the ($n\times n)$ matrix $\Pi$ has rank $r$, $0<r<n$, it has $n-r$ non null eigenvalues, implying that the series in ${\bf Y}_t$ has at least one unit root and its valid representation is given by the VECM in Equation \eqref{eq:vec}.
In this case, $\Pi=\alpha\beta'$, where $\alpha$ and $\beta$ are matrices ($n\times r$) of rank $r$.
Matrix $\beta$ denotes the one with the cointegrating vectors and matrix $\alpha$ is called the loading matrix, since it contains the weights of the equilibrium relationships.
The tests developed in \cite{johansen1990} focus on the rank of matrix $\Pi$.

The pioneer Bayesian works to study VAR models and reduced rank regressions are \cite{dejong1992, bauwens1996,geweke1996}. 
However, the main concern of these papers is to estimate the model parameters and their (marginal) posterior distributions. 
The usual approach is to assume a given rank for the long run matrix $\Pi$, and proceed with all the computations conditional on the given rank. 
The Bayesian initiatives to test the rank of the referred matrix are recent, the main reference for Bayesian inference on VECM's being \cite{koop2006}.

To justify inferential procedures based on prespecified ranks of matrix $\Pi$, \cite{bauwens1999} argued that an empirical cointegration analysis should be based on economic theory, which proposes models obeying equilibrium relationships.
According to this view, cointegration research should be ``confirmatory'' rather than ``exploratory''.
Even though the advocated conditional inference is of simple implementation and very useful for small samples, \cite{bauwens1999} recognized that tests for the rank of matrix $\Pi$ should be developed. 
To our knowledge, few initiatives with this purpose were developed up to now. 

One common approach to test sharp hypotheses in the Bayesian framework is by means of Bayes factors. 
Testing the rank of matrix $\Pi$ by Bayes factors implies several computational complications and requires the use of proper priors, as shown in \cite{kleibergen1997}. 
Following an informal approach, \cite{bauwens1996} obtained the posterior distribution of the ordered eigenvalues of the ``squared''  long run matrix, $\Pi'\cdot\Pi$, obtained from a VAR model without assuming the existence of cointegration relations. 
As the long run matrix has a reduced rank, it has some null eigenvalues, and this should be revealed by the fact that the smallest eigenvalues should have a lot of probability mass accumulated on values close to zero. 
The computations can be made straightforwardly, simulating values for the long run matrix from its (marginal) posterior distribution, which is a matrix $t$-Student distribution under the non informative prior \eqref{eq:prior_var}, also considered in the sequel.

Another common procedure is to estimate the rank of $\Pi$ as the value $r$ that maximizes the (marginal) posterior distribution of the rank.
Conditioned on such an estimate, one proceeds to derive the full posterior and eventually estimate the cointegration space, i.e., the linear space spanned by $\beta$.

A different approach was proposed by \cite{chao1999}, who used the Posterior Information Criterion (PIC), developed in \cite{phillips1996}, as a criterion to choose the mode of the posterior distribution of the rank of $\Pi$. 
However, as highlighted in \cite{koop2006}, one of the advantages of the Bayesian approach is the possibility to incorporate the uncertainty about the parameters in the analysis, represented by the posterior distribution of the rank and, whatever the tool the scientist uses to infer the value of $r$, it is derived from this posterior distribution.

The authors of \cite{kleibergen2002} nested the reduced rank models in an unrestricted VAR and used Metropolis--Hastings sampling with the Savage--Dickey density ratio---see \cite{verdinelli1995}---to estimate the Bayes Factors of all the models with incomplete rank up to the model with full rank. 
The Bayes Factor derivation requires the estimation of an error correction factor for the incomplete rank.  
This factor, however, is not defined for improper priors due to a problem known as \textit{Bartlett paradox}, which arises whenever one compares models of different dimensions.  
The difficulty is relevant in the present case because, after deriving the rank posterior density, one may consider that models of different dimensions are being compared. 
The paradox is stated informally as: improper priors should be avoided when one computes Bayes Factors (except for parameters common to both models) as they depend on arbitrary constants (that are integrals).

More recently, \cite{villani2005} developed an efficient procedure to obtain the posterior distribution of the rank using a uniform proper prior over the cointegration space linearly normalized. 
The author derived solutions for the posterior probabilities for the null rank and for the full rank of $\Pi$. 
The posterior probabilities of each intermediate rank are derived from the posterior samples of the matrices that compose the long run matrix ($\alpha$ and $\beta)$, properly normalized, under each rank and using the method proposed by \cite{chib1995}.

\section{Implementing the FBST as a Cointegration Test}
\label{sec:impl_coint}

This section describes how to implement the FBST to test for cointegration. 
We will proceed in the same spirit of Section \ref{sec:FBST_unit_root}, i.e., describing the  steps given in Table \ref{tab:pseudocode} to implement the test for cointegration.

Let us begin recalling the VECM given by Equation \eqref{eq:vec}:

\begin{equation}
\Delta {\bf Y}_t={\bf c}+\Phi_0 {\bf D}_t+\Gamma_1 \Delta {\bf Y}_{t-1}+\ldots+\Gamma_{p-1}\Delta {\bf Y}_{t-p+1}+\Pi {\bf Y}_{t-1}+{\bf E}_t, \tag{\ref{eq:vec}}
\end{equation}

\noindent
$t=1,\ldots,T+p$, in which ${\bf E}_t\stackrel{i.i.d.}\sim N_n({\bf 0},\Sigma)$ with ${\bf 0}$ a null  vector of dimension $n\times 1$ and $\Omega$ a symmetric positive definite real matrix.
Notice that these assumptions already specify the statistical model (Gaussian) and its implied likelihood.
Before giving it explicitly, let us rewrite Equation \eqref{eq:vec} using matrix notation: 



\begin{equation}
\Delta {\bf Y}={\bf Z}\cdot\eta+{\bf E} 
\label{eq:mat}
\end{equation}

\noindent
where $\Delta {\bf Y}=\left[ \begin{array}{c}\Delta {\bf Y}_{p+1}' \\ \Delta {\bf Y}_{p+2}'\\ \vdots \\ \Delta {\bf Y}_{T+p}' \end{array} \right]$, ${\bf Z}=\left[ \begin{array}{cccccc}
1 &{\bf D}_{p+1}' &\Delta {\bf Y}_{p}'&\ldots&\Delta {\bf Y}_{2}'&{\bf Y}_p'\\
1 & {\bf D}_{p+2}'& \Delta {\bf Y}_{p+1}'&\ldots&\Delta {\bf Y}_{3}'&{\bf Y}_{p+1}'\\
\vdots&\vdots&\vdots&&\vdots&\vdots\\
1&{\bf D}_{T+p}'&\Delta {\bf Y}_{T-1}'&\ldots&\Delta {\bf Y}_{T+p-1}'&{\bf Y}_{T+p-1}'\\
\end{array} \right]$, $\eta=\left[ \begin{array}{c}{\bf c}'\\ \Phi_0\\ \Gamma_1 \\ \vdots \\ \Gamma_{p-1} \\ \Pi\end{array} \right]$ and the error vector is given by ${\bf E} \sim MN_{T\times n}(0, I_T, \Omega)$, denoting the matrix normal distribution.
See Appendix \ref{app:dist} for more information on this distribution.
Now the parameter vector is given by $\Theta=(\eta, \Omega)$.

Notice that $\Delta {\bf Y}$ is formed by piling up $T$ transposed vectors $\Delta {\bf Y}_t$, thus resulting in a matrix with $T$ lines and $n$ columns ($n$ is the number of time series in vector ${\bf Y}_t$), those being also dimensions of matrix ${\bf E}$.
Matriz ${\bf Z}$ is constructed likewise---always piling up the transposed vectors---resulting in a matrix with $T$ lines and $pn+n+1$ columns.
Finally, matrix $\eta$ has the matrices of coefficients, all piled up properly, resulting in a matrix with $pn+n+1$ lines and $n$ columns.

Given the assumptions above, $\Delta{\bf Y}\sim MN_{T\times n}({\bf Z}\cdot\eta, I_T, \Omega)$, implying that the likelihood is 

\begin{displaymath}
L(\Theta\mid{\bf y})\propto |\Omega|^{-T/2}\textrm{exp}\left\{-\frac{1}{2}\cdot \textrm{tr}[\Omega^{-1}(\Delta {\bf Y}-{\bf Z}\cdot \eta)'(\Delta {\bf Y}-{\bf Z}\cdot \eta)]\right\}
\end{displaymath}

\noindent
where ${\bf y}$ denotes the set of observed values of vectors ${\bf Y}_t$ for $t=1,\ldots,T+p$.
As in Section \ref{sec:FBST_unit_root}, we will consider an improper prior for $\Theta$, given by

\begin{equation}\label{eq:prior_var}
h(\Theta)=h(\eta,\Omega)\propto |\Omega|^{-(n+1)/2},
\end{equation}

\noindent
and our reference density, $r(\Theta)$, will be proportional to a constant, leading to a surprise function equivalent to the (full) posterior distribution.
These choices correspond to steps 1 and 2 of Table \ref{tab:pseudocode}.
These modeling choices imply the following posterior density:

\begin{equation}
\begin{array}{ccl}
g(\Theta\mid{\bf y})&\propto& |\Omega|^{-(T+n+1)/2}\textrm{exp}\left\{-\frac{1}{2}\cdot \textrm{tr}[\Omega^{-1}(\Delta {\bf Y}-{\bf Z}\cdot \eta)'(\Delta {\bf Y}-{\bf Z}\cdot \eta)]\right\}\\
&=&|\Omega|^{-(T+n+1)/2}\textrm{exp}\left\{-\frac{1}{2}\cdot \textrm{tr}\{\Omega^{-1}[{\bf S}+(\eta-\widehat{\eta})'\cdot{\bf Z}'{\bf Z}\cdot(\eta-\widehat{\eta})]\}\right\}
\end{array} \label{eq:pos}
\end{equation}

\noindent
where $\widehat{\eta}=({\bf Z}'{\bf Z})^{-1}{\bf Z}'\Delta {\bf Y}$ and ${\bf S}=\Delta {\bf Y}'\Delta {\bf Y}-\Delta {\bf Y}'{\bf Z}({\bf Z}'{\bf Z})^{-1}{\bf Z}'\Delta {\bf Y}$.

To implement Step 3 of Table \ref{tab:pseudocode}, we need to find the maximum a posteriori of \eqref{eq:pos} under the constraint $\Theta\subset\Theta_H$, i.e., we need to maximize the posterior in $\Theta_H$.
Since we are testing the rank of matrix $\Pi$, as discussed in the beginning of Section \ref{sec:bayes_coint}, it is necessary to maximize the posterior assuming the rank of $\Pi$ is $r$, $0\leq r\leq n$.
Thanks to the modeling choices made here---Gaussian likelihood and Equation \eqref{eq:prior_var} as prior---our posterior is almost identical to a Gaussian likelihood, allowing us to find this maximum using a strategy similar to that proposed by \cite{johansen1990}, who derived the maximum of the (Gaussian) likelihood function  assuming a reduced rank for $\Pi$.
We will summarize Johansen's algorithm, providing in Appendix \ref{app:johansen} a heuristic argument of why it indeed provides the maximum value of the posterior under the assumed hypotheses.

We begin estimating a VAR($p-1$) model for $\Delta {\bf Y}_t$ with all the explanatory variables shown in \eqref{eq:vec} except for ${\bf Y}_{t-1}$.
Using the matrix notation established above, this corresponds to estimate

\[
\Delta {\bf Y}={\bf Z}_1\cdot \eta_1+{\bf U},
\]

\noindent
where ${\bf Z}_1=\left[ \begin{array}{ccccc}
1 &{\bf D}_{p+1}' &\Delta {\bf Y}_{p}'&\ldots&\Delta {\bf Y}_{2}'\\
1 & {\bf D}_{p+2}'& \Delta {\bf Y}_{p+1}'&\ldots&\Delta {\bf Y}_{3}'\\
\vdots&\vdots&\vdots&&\vdots\\
1&{\bf D}_{T+p}'&\Delta {\bf Y}_{T-1}'&\ldots&\Delta {\bf Y}_{T+p-1}'\\
\end{array} \right]$ and $\eta_1=\left[ \begin{array}{c}{\bf e}'\\ \tau_0\\ \upsilon_1 \\ \vdots \\  \upsilon_{p-1} \end{array} \right]$ showing that ${\bf Z}_1$ is obtained from matrix ${\bf Z}$ extracting its last $n$ columns, exactly those corresponding to ${\bf Y}_{t-1}$.

We also estimate a second set of auxiliary equations, regressing ${\bf Y}_{t-1}$ on a vector of constants and ${\bf D}_t$, $\Delta {\bf Y}_{t-1}$, \ldots, $\Delta {\bf Y}_{t-p+1}$. 
By piling up all the (transposed) vectors ${\bf Y}_{t-1}'$ for $t=p+1,\ldots, T+p$, we have a $(T\times n)$ matrix, denoted by ${\bf Y}_{-1}$.
As above, these equations can be represented by

\[
{\bf Y}_{-1}={\bf Z}_1\cdot \eta_2+{\bf V},
\]

\noindent
where ${\bf Y}_{-1}=\left[ \begin{array}{c}{\bf Y}_{p}' \\ {\bf Y}_{p+1}'\\ \vdots \\ {\bf Y}_{T+p-1}' \end{array} \right]$ and $\eta_2=\left[ \begin{array}{c}{\bf m}'\\ \nu_0\\ \zeta_1 \\ \vdots \\  \zeta_{p-1} \end{array} \right]$.

Considering the OLS estimates of these sets of equations and their respective estimated residuals, we may write

\begin{align}
\widehat{\Delta {\bf Y}}&={\bf Z}_1\cdot \widehat{\eta}_1+\widehat{{\bf U}}\\
\widehat{{\bf Y}}_{-1}&={\bf Z}_1\cdot \widehat{\eta}_2+\widehat{{\bf V}}
\end{align}

\noindent
where $\widehat{\eta}_1=({\bf Z}_1'{\bf Z}_1)^{-1}{\bf Z}_1'\cdot \Delta {\bf Y}$, ~  $\widehat{\eta}_2=({\bf Z}_1'{\bf Z}_1)^{-1}{\bf Z}_1'\cdot {\bf Y}_{-1}$, ~ $\widehat{{\bf U}}$ and $\widehat{{\bf V}}$ are the respective matrices of estimated residuals. 
Thanks to the Frisch--Waugh--Lovell theorem---see \cite{greene2008} theorem 3.3 or \cite{davidson2004} Section 2.4---, it is possible to show that the estimated residuals of these auxiliary regressions are related by $\Pi$ in the following regressions:

\begin{equation}\label{eq:fwl}
\widehat{{\bf U}}=\Pi ~ \widehat{{\bf V}}+\widehat{{\bf W}}.
\end{equation}

One can prove that the OLS estimates of $\Pi$ obtained from \eqref{eq:mat} and from \eqref{eq:fwl} are numerically identical, as the estimated residuals $\widehat{{\bf E}}$ and $\widehat{{\bf W}}$.

The second stage of Johansen's algorithm requires the computation of the following sample covariance matrices of the OLS residuals obtained above:

\begin{align*}
\widehat{\Sigma}_{{\bf V}{\bf V}}={1\over T}\cdot \widehat{{\bf V}}'\widehat{{\bf V}} & \qquad \widehat{\Sigma}_{{\bf U}{\bf U}}={1\over T}\cdot \widehat{{\bf U}}'\widehat{{\bf U}} \\   
\widehat{\Sigma}_{{\bf U}{\bf V}}={1\over T}\cdot \widehat{{\bf U}}'\widehat{{\bf V}} & \qquad \widehat{\Sigma}_{{\bf V}{\bf U}}=\widehat{\Sigma}_{{\bf U}{\bf V}}'
\end{align*}

\noindent
and, from these, we find the $n$ eigenvalues of matrix

\[
\widehat{\Sigma}_{{\bf V}{\bf V}}^{-1}\cdot \widehat{\Sigma}_{{\bf V}{\bf U}}\cdot\widehat{\Sigma}_{{\bf U}{\bf U}}^{-1}\cdot\widehat{\Sigma}_{{\bf U}{\bf V}},
\]

\noindent
ordering them decreasingly $\widehat{\lambda}_1>\widehat{\lambda}_2>\ldots>\widehat{\lambda}_n$.
The maximum value attained by the log posterior subject to the constraint that there are $r$ ($0\leq r\leq n$) cointegration relationships is

\begin{equation}\label{eq:pos_r}
\ell^*=K-{(T+n+1)\over 2}\cdot\log|\widehat{\Sigma}_{{\bf U}{\bf U}}|-{T+n+1\over 2}\cdot\sum_{i=1}^r\log(1-\widehat{\lambda}_i),
\end{equation}

\noindent
where $K$ is a constant that depends only on $T$, $n$ and ${\bf y}$ by means of the marginal distribution of the data set, ${\bf y}$.
Since $\ell^*$ represents the maximum of the log-posterior, to obtain $s^*$, one should take $s^*=\exp(\ell^*)$, completing step 3 of Table \ref{tab:pseudocode}.

As in Section \ref{sec:FBST_unit_root}, we compute $\overline{ev}$ in step 4 by means of a Monte Carlo algorithm. 
It is easy to factor the full posterior \eqref{eq:pos} as a product of a (matrix) normal and an Inverse-Wishart, suggesting a Gibbs sampler to generate random samples from the full posterior.
See Appendix \ref{app:dist} for more on the Inverse-Wishart distribution.
Thus, the conditional posteriors for $\eta$ and $\Omega$ are, respectively, 

\begin{align}
g_{\eta}(\eta\mid\Omega,{\bf y})&\propto MN_{n\times k}(\widehat{\eta}, ({\bf Z}'{\bf Z})^{-1},\Omega)\label{eq:cond_eta}\\
g_{\Omega}(\Omega\mid\eta,{\bf y})&\propto IW(\Omega|{\bf S}+(\eta-\widehat{\eta})'\cdot{\bf Z}'{\bf Z}\cdot(\eta-\widehat{\eta}), T)\label{eq:cond_Omega}
\end{align}

\noindent
where ${\bf S}=\Delta {\bf Y}'\Delta {\bf Y}-\Delta {\bf Y}'{\bf Z}({\bf Z}'{\bf Z})^{-1}{\bf Z}'\Delta {\bf Y}$, $IW$ denotes the Inverse-Wishart, $k=pn+n+1$ is the number of lines of $\eta$, and $\widehat{\eta}$ its OLS estimator, as above.
From a Gibbs sampler set with these conditionals, we obtain a random sample from the full posterior to estimate $\overline{ev}$ as the proportion of sampled vectors that generate a value for the posterior greater than $s^*$.
Finally, we obtain $ev=1-\overline{ev}$ in the final step (5).
The whole implementation for cointegration tests, following the assumptions made in this section, is summarized in Table \ref{tab:coint_tests}.
See Appendix \ref{app:computer} for more information on the computational resources needed to implement the steps given by Table \ref{tab:coint_tests}.

\begin{table}[!h]
\centering
\caption{Pseudocode to implement the FBST to cointegration tests}
\begin{tabular}{l}
\toprule
\textbf{General algorithm}: compute $ev$ supporting hypothesis $H: \text{rank}(\Pi)=r$ ($0\leq r\leq n)$ in model \eqref{eq:vec}\\
\midrule
1. Statistical model: Gaussian; prior: $h(\Theta)\propto|\Omega|^{-(n+1)/2}$.\\
2. Reference density: $r(\Theta)\propto 1$; relative surprise function: $g(\Theta\mid{\bf y})$.\\
3. Find $s^*$: Johansen's algorithm; obtain $\ell^*$ from Equation \eqref{eq:pos_r} with $s^*=\exp(\ell^*)$.\\
4. Gibbs sampler (from Equations \eqref{eq:cond_eta} and \eqref{eq:cond_Omega}) to obtain $N$ random samples of parameter vectors from \eqref{eq:pos}.\\
~ ~ Evaluate the posterior at the sampled vectors and estimate $\overline{ev}$ as the proportion of $N$ for which the \\
~ ~ evaluated values are larger than $s^*$.\\
5. Find $ev=1-\overline{ev}$.\\
\hline
\end{tabular}
\label{tab:coint_tests}
\end{table}

Before presenting the results of the procedure applied to real data sets, it is important to remark one feature of the FBST applied to cointegration tests.
The estimated eigenvalues of matrix $\Pi$, $\widehat{\lambda}_i$, correspond to the squared canonical correlations between $\Delta {\bf Y}_t$ and ${\bf Y}_{-1}$ corrected for the variable in ${\bf Z}_1$ and therefore lie between 0 and 1.
Therefore, \eqref{eq:pos_r} shows that $\ell^*_{0}\leq\ell^*_{1}\leq \ldots \ell^*_{n}$, where $\ell^*_r$ denotes the maximum of the posterior \eqref{eq:vec} assuming $\Pi$ has rank $0\leq r\leq n$. 
Therefore, one may say that the hypotheses
rank($\Pi)=r$ are nested, in the sense that the respective e-values obtained by the FBST for these hypotheses are always non-decreasing $ev(0)\leq ev(1)\leq \ldots \leq ev(n)$.

This nested formulation is also present in the frequentist procedure proposed by \cite{johansen1990}, based on the likelihood ratio statistics for successive ranks of $\Pi$.
Thus, the FBST should be used, like the maximum eigenvalue test, in a sequential procedure to test for the number of cointegrating relationships.
We will show how this should be done in presenting the applied results in the sequel.  

\subsection{Results}

Now we present, by means of four examples, the application of FBST as a cointegration test.
In all the examples, we have adopted a Gaussian likelihood and the improper prior \eqref{eq:prior_var}. 
The Gibbs sampler was implemented as described above, providing 51,000 random vectors from the posterior \eqref{eq:pos}.
The first 1000 samples were discarded as a burn-in sample, the remaining 50,000 being used to estimate the integral \eqref{eq:evc}.
The tables show the e-value computed from the FBST and the maximum eigenvalue test statistics with their respective $p$-values.




\

\noindent
\textbf{Example 1}. \emph{We analyzed four electroencephalography (EEG) signals from a subject that has previously presented epileptic seizures.
The original study, \cite{shoeb2009}, had the aim of detecting seizures based on multiple hours of recordings for each individual and the cointegration analysis of the mentioned signals was presented by \cite{ostergaard2017}.
In fact, the cointegration hypothesis is tested using the phase processes estimated from the original signals.
This is done by passing the signal into the Hilbert transform and then ``unwrapping'' the resulting transform. 
Sections 2 and 5 of \cite{ostergaard2017,freeman2007} provide more details on the Hilbert transform and unwrapping.}\\

\begin{figure}[!h]
\centering
\includegraphics[width=.8\linewidth]{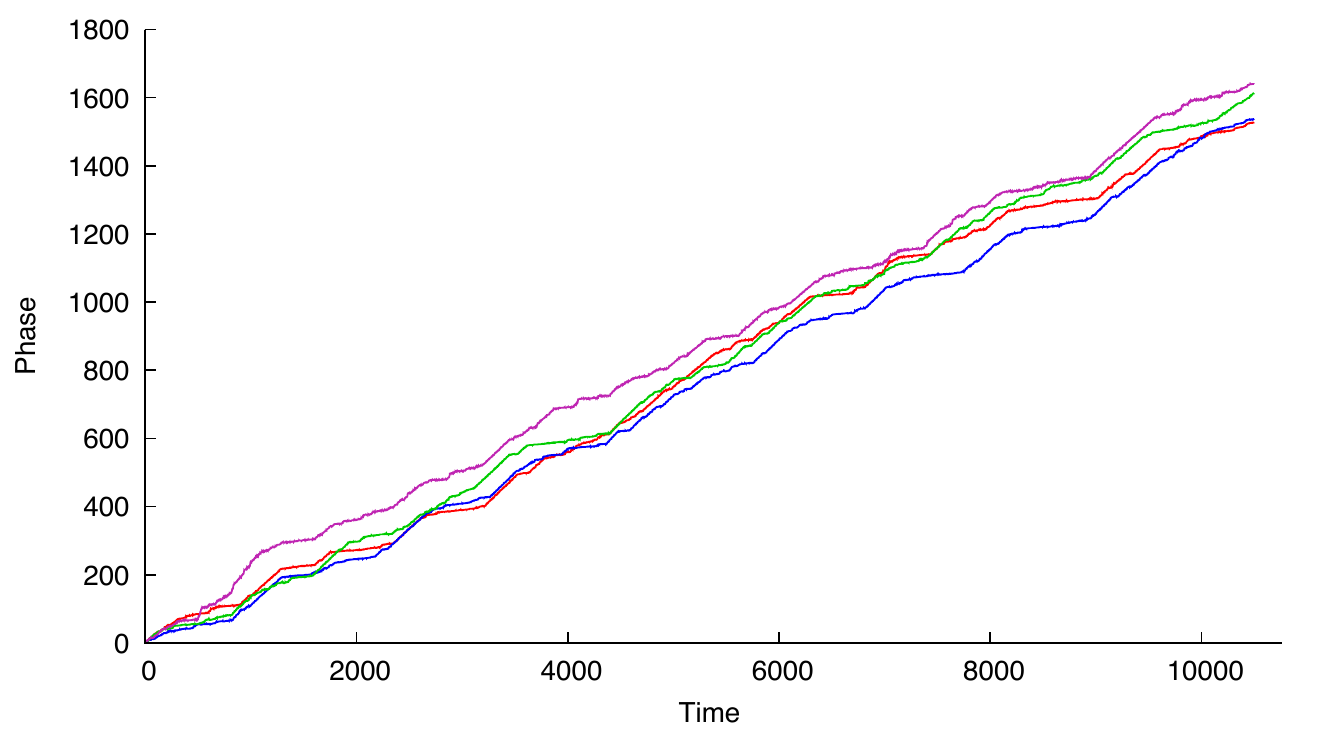}
\caption{Estimated phase processes prior to a seizure.}
\label{fig:prior}

\end{figure}

\begin{figure}[!h]
\centering
\includegraphics[width=.8\linewidth]{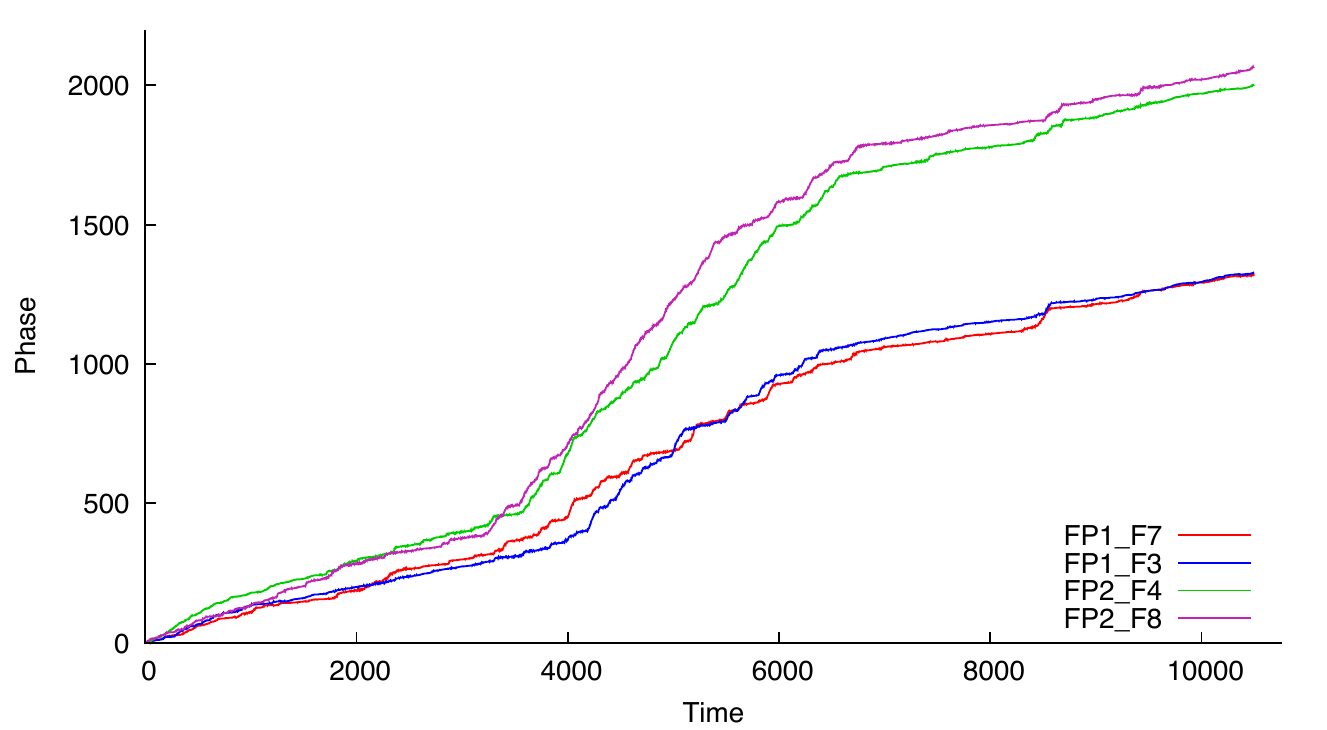}
\caption{Estimated phase processes during a seizure.}
\label{fig:during}

\end{figure}

The labels of the modeled series refer to the electrodes on the scalp. As seen in \mbox{Figures \ref{fig:prior} and \ref{fig:during}}, the series are called FP1-F7, FP1-F3, FP2-F4, and FP2-F8, where FP refers to the frontal lobes and F refers to a row of electrodes placed behind these. 
Even numbered electrodes are on the right side and odd numbered electrodes are on the left side. 
The electrodes for these four signals mirror each other on the left and right sides of the scalp. 
The recordings of the studied subject, an 11-year-old female, identified a seizure in the interval (measured in seconds) $[2956, 2996]$. 
Therefore, like \cite{ostergaard2017}, we analyze the period of 41 seconds prior to the seizure---interval $[2956, 2996]$---and the subsequent 41 seconds---interval $[2996, 3036]$---the seizure period.
In the sequel, we will refer to these as {\it prior to seizure} and {\it during seizure}, respectively.
Since the sample frequency has 256 measurements per second, there are a total of 10,496 measurements for each of the four signals.
\cite{ostergaard2017} used 40 seconds for each period, obtaining slightly different results.

Figures \ref{fig:prior} and \ref{fig:during} display the estimated phases based on the original signals.
The model proposed by \cite{ostergaard2017} is a VAR(1), resulting in a VECM given by

\begin{equation}\label{eq:ex_eeg}
\Delta{\bf Y}_t={\bf c}+\Pi{\bf Y}_{t-1}+{\bf E}_t.
\end{equation}

Tables \ref{tab:prior} and \ref{tab:during} present the results that essentially lead to the same conclusions obtained by \cite{ostergaard2017}, even though they have based their findings on the trace test.
See Table 8 of \cite{ostergaard2017}.
\unskip
\unskip
\begin{table}[!h]\centering
{

\caption{FBST and max. eig. test: prior to seizure}\label{tab:prior}
\begin{tabular}{lccc}
\toprule
\boldmath{$H_0$}	&\textbf{FBST}  	    & \textbf{Max.}	&\boldmath{$p$}\textbf{-Value}\\ \midrule
$r=0$	& $\simeq$$0$	& 60.966 & $\simeq$$0$ \\
$r=1$	& 0.0691		& 30.727 & 0.0010 \\
$r=2$	& 0.9990		& 11.458 & 0.1337 \\
$r=3$	& $\simeq$$1$	& 0.0812 & 0.7757 \\
\hline
\end{tabular}
}\end{table}
\unskip
\begin{table}[!h]\centering
\qquad
{

\caption{FBST and max. eig. test: during seizure}\label{tab:during}
\begin{tabular}{lccc}
\toprule
\boldmath{$H_0$}	&\textbf{FBST}  	    & \textbf{Max.}	&\boldmath{$p$}\textbf{-Value}\\ \midrule
$r=0$	& $\simeq$$0$	& 1120.5 & $\simeq$$0$ \\
$r=1$	& 0.1144		& 31.563 & 0.0007 \\
$r=2$	& 0.9999		& 6.5015 & 0.5574 \\
$r=3$	& $\simeq$$1$	& 1.4383 & 0.2304 \\
\hline
\end{tabular}
}
\end{table}

The comparison between $p$-values and the FBST e-values must be made carefully, the main reason being the fact that $p$-values are not measures supporting the null hypothesis, while e-values provide exactly such a kind of support. 
That being said, a possible way to compare them is by checking the decision their use recommend regarding the hypothesis being tested, i.e., to reject or not the null hypothesis.

Frequentist tests often adopt a significance level approach: given an observed $p$-value, the hypothesis is rejected if the $p$-value is smaller or equal to the mentioned level, usually 0.1, 0.05, or 0.01.
Since the cointegration ranks generate nested likelihoods, the hypotheses are tested sequentially, starting with null rank, $r=0$.
For Table \ref{tab:prior}, adopting a 0.01 significance level, the maximum eigenvalue test would reject $r=0$ and $r=1$, and would not reject $r=2$.
The same conclusions follow for Table~\ref{tab:during}.
Thus, the recommended action is to work, for estimation purposes for instance, assuming two cointegration relationships.

The question on which threshold value to adopt for the FBST was already mentioned on Section~\ref{subsec:practical}, but it is worthwhile to underline it once more.
We highly recommend a principled approach deriving the cut-off value from a loss function, which is specific for the problem at hand and the purposes of the analysis.
A naive but simpler approach would be to reject the hypothesis if the e-value is smaller than 0.05 or 0.01, emulating the frequentist strategy.
Even not recommending this path, since $p$-values are not supporting measures for the hypothesis being tested while e-values are, the researcher may numerically compare $p$-values and e-values in a specific scenario.
If the researcher derived the $p$-values from a generalized likelihood ratio test, it is possible to asymptotically compare them.
The relationship is: $ev=1-F_m[F^{-1}_{m-h}(1-{\bf p})]$, where $m$ is the dimension of the full parameter space, $h$ the dimension of the parameter space under the null hypothesis, $F_m$ the chi-square distribution function with $m$ degrees of freedom and {\bf p} the corresponding $p$-value. See \cite{diniz2012b, stern2013} for the proof of the asymptotic relationship between e-values and $p$-values.

Since the maximum eigenvalue test is derived as a likelihood ratio test, this comparison may be done for the results of all the examples presented here, and more appropriately to this example, given its sample size of 10,496 observations.
Regarding Tables \ref{tab:prior} and \ref{tab:during}, one could be in doubt regarding whether to reject or not the hypothesis $r=1$ since the e-values are larger than 0.01.
However, for this model and hypothesis, the e-value corresponding to 0.01 is 0.436.
Therefore, in both tables, one could reject the hypothesis and proceed to the next rank that has plenty of evidence in its favor.
In conclusion, the practical decisions of both tests (FBST and maximum eigenvalue) would be the same: to not reject $r=2$.

\

\noindent
\textbf{Example 2} (\cite{turasie2016}). \emph{Compare three methods for modeling empirical seasonal temperature forecasts over South America. 
One of these methods is based on a (possible) long-term cointegration relationship between the temperatures of the quarter March--April--May (MAM) of each year and the temperature of the previous months of November--December--January (NDJ).
When there is such a relationship, the authors used the NDJ temperatures (of the previous year) as a predictor for the following MAM season.}\\

The original data set has monthly temperatures for each coordinate (latitude and longitude) of the covered area.
The mentioned series of temperatures (MAM and NDJ) are computed as seasonal averages from this monthly data set by averaging over consecutive three months. 
Since we have data available from January 1949 to May 2020, the time series of monthly and seasonal average surface temperatures of length 72 for each grid point.

The authors of \cite{turasie2016} consider ${\bf Y}_t$ as a two-dimensional vector, its first component being the seasonal (average) MAM temperature of year $t$ and the second component the seasonal NDJ temperature of the {\it previous} year.
They consider a VAR(2) without deterministic terms to model the series, resulting in a VECM

\begin{equation}
\Delta {\bf Y}_t=\Gamma_1\Delta{\bf Y}_{t-1}+\Pi{\bf Y}_{t-1}+{\bf E}_t.
\end{equation}

\begin{figure}[!h]
\begin{minipage}[b]{0.48\linewidth}
\includegraphics[width=\linewidth]{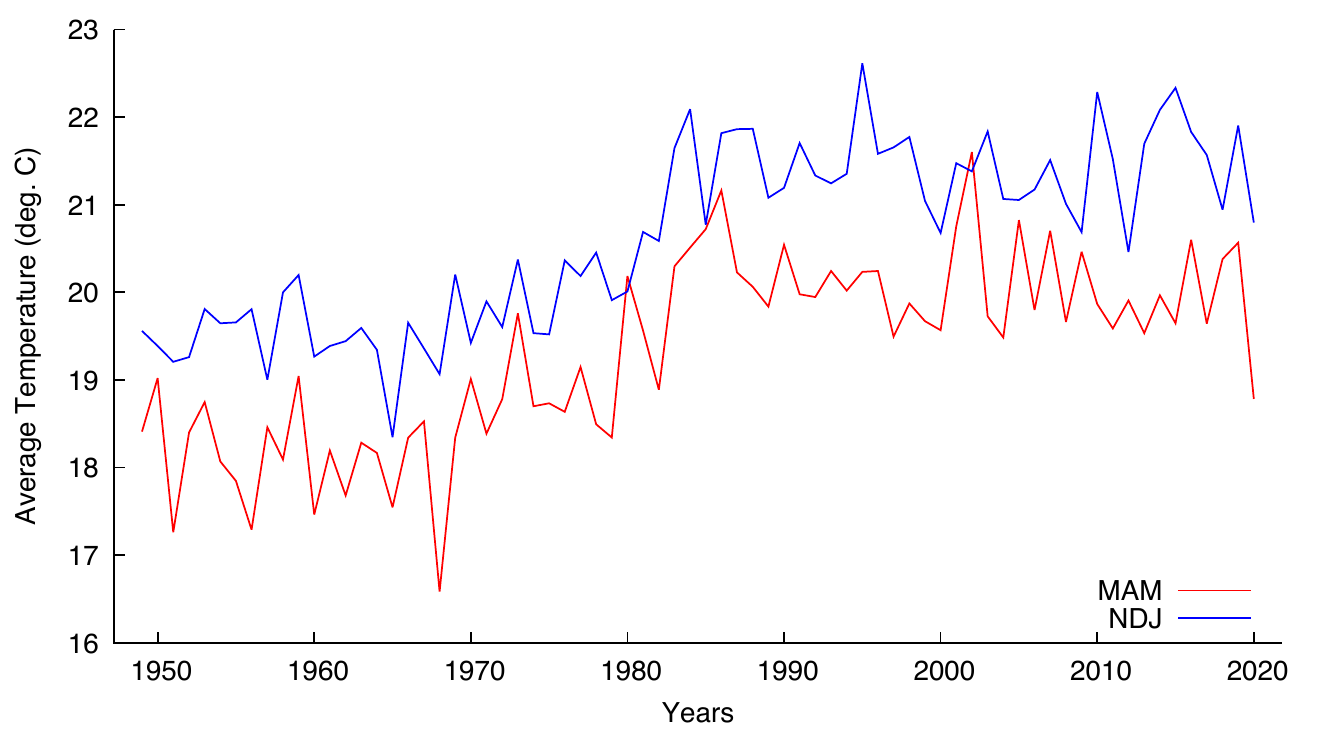}
\caption{Seasonal (MAM and NDJ) temperatures for S\~ao Paulo from 1949 to 2020}
\label{fig:sp}
\end{minipage} \hfill
\begin{minipage}[b]{0.48\linewidth}
\includegraphics[width=\linewidth]{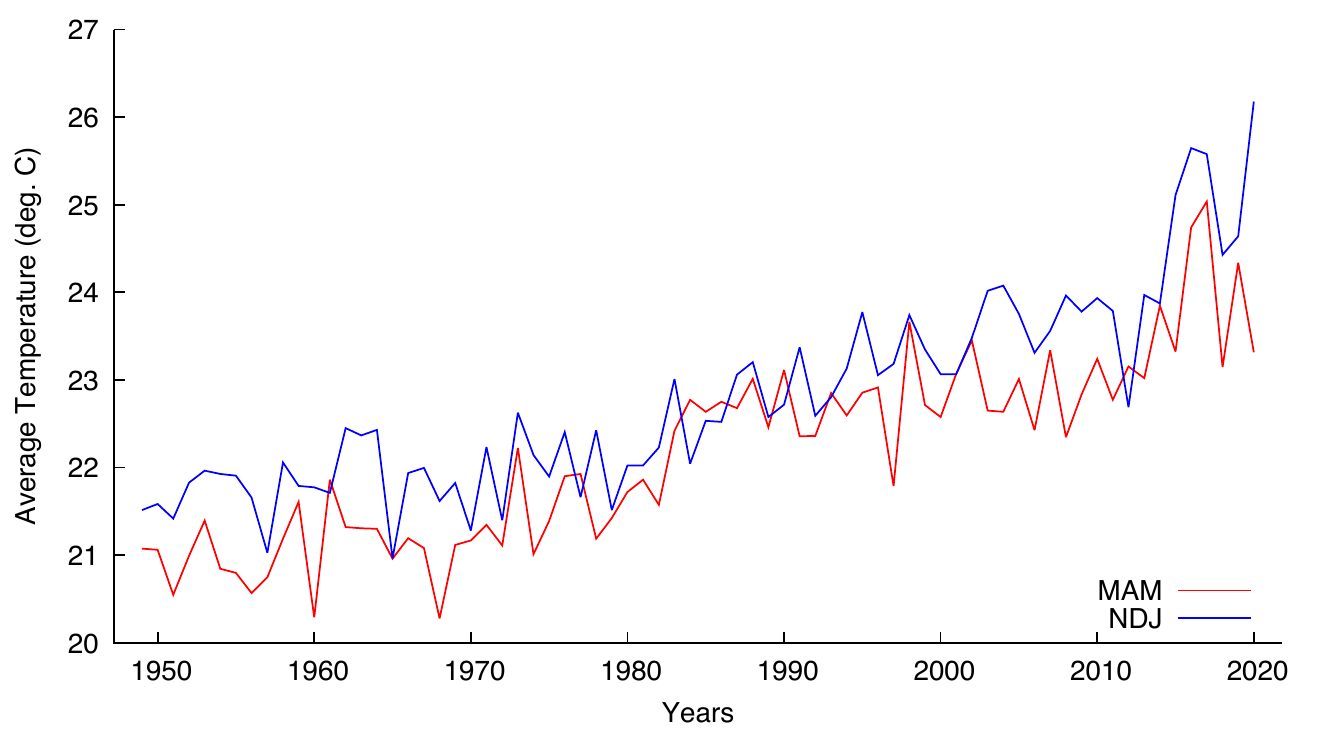}
\caption{Seasonal (MAM and NDJ) temperatures for Bras\'ilia from 1949 to 2020}
\label{fig:br}
\end{minipage}
\end{figure}

We have chosen five grid points corresponding to major Brazilian cities to test the cointegration hypothesis of the mentioned seasonal series.
The coordinates chosen were the closest ones from: 23.5505$^{\circ}$ S, 46.6333$^{\circ}$ W for S\~ao Paulo; 22.9068$^{\circ}$ S, 43.1729$^{\circ}$ W for Rio de Janeiro; 19.9167$^{\circ}$ S, 43.9345$^{\circ}$ W for Belo Horizonte; 15.8267$^{\circ}$ S, 47.9218$^{\circ}$ W for Bras\'ilia and 12.9777$^{\circ}$ S, 38.5016$^{\circ}$ W for Salvador. 
Figures \ref{fig:sp} and \ref{fig:br} show the seasonal temperatures for S\~ao Paulo and Bras\'ilia, respectively, indicating that the cointegration hypothesis is plausible for both cities.

\begin{table}[!h]
\centering
\caption{FBST and maximum eigenvalue test applied to temperature data (MAM and NDJ series) of the mentioned Brazilian cities}
\begin{tabular}{lccccccc}
\hline
                    &\multicolumn{3}{c}{$H_0:r=0$}& &\multicolumn{3}{c}{$H_0:r=1$}\\ \cline{2-4} \cline{6-8}
Cities		        &FBST  	    & Max.	&$p$-value&  &FBST  	    & Max.	&$p$-value\\ \hline
S\~ao Paulo	        &0.0012 &  33.302  & $\simeq$$0$ &&$\simeq$$1$ & 0.0893 & 0.8205	\\
Rio de Janeiro		&0.0273 & 23.294  &	0.0004 && $\simeq$$1$ & 2.43e-5 & 0.9986\\ 
Belo Horizonte      &0.0173 & 24.621  & 0.0001 && $\simeq$$1$ & 0.0963 & 0.8126\\
Bras\'ilia          &0.1129& 18.008  & 0.0045 && 0.9999 & 1.3321 & 0.2892 \\
Salvador            &0.0172 & 24.431  & 0.0001 && $\simeq$$1$ & 0.2450 & 0.6838 \\ \hline
\end{tabular}
\label{tab:ex_clim}
\end{table}

The results are shown in Table \ref{tab:ex_clim}.
Assuming a significance level of 0.01, the maximum eigenvalue test reject the null rank and do not reject $r=1$ for all the five cities.
If we adopt the asymptotic relationship between $p$-values and e-values for the model under analysis, we obtain an e-value of 0.276 corresponding to a 0.01 $p$-value for $r=0$.
Therefore, the FBST would also reject the null rank for all the cities.
The hypothesis $r=1$ is not rejected since all the e-values are close to 1, once more agreeing with the maximum eigenvalue test.

One remark about Bras\'ilia seems in order.
The city was built to be the federal capital, being officially inaugurated on 21 April 1960.
The construction began circa 1957 and before that the site had no human occupation.
The process of moving all the administration from Rio de Janeiro, the former capital, was slow and only the 1980 census detected a population over 1 million inhabitants.
The present population is almost 3,2 million inhabitants living in the Federal District that includes Bras\'ilia and minor surrounding cities.
Figure \ref{fig:br} indicates that the seasonal temperatures began to rise exactly after 1980.

\

\noindent{\textbf{Example 3}: \emph{we applied the FBST to the Finish data set used in their seminal work} \cite{johansen1990}.}\\

The authors used the logarithms of the series of $M1$ monetary aggregate, inflation rate, real income, and the primary interest rate set by the Bank of Finland to model the money demand, which, in theory, follows a long-term relationship. 
The sample has 106 quarterly observations of the mentioned variables, starting at the second trimester of 1958 and finishing in the third trimester of 1984.
The chosen model was a VAR(2) with unrestricted constant, meaning that the series in ${\bf Y}_t$ have one unit root with drift vector ${\bf c}$ and the cointegrating relations may have a non-zero mean. For more information about how to specify deterministic terms in a VAR see \cite{helmut2005}, chapter 6. 
Seasonal dummies for the first three quarters of the year were also considered in the model chosen by \cite{johansen1990}.
Writing the model in the error correction form, we have:

\begin{equation}\label{eq:finland}
\Delta {\bf Y}_t={\bf c}+\Phi_{0,1}{\bf D}_{1t}+\Phi_{0,2}{\bf D}_{2t}+\Phi_{0,3}{\bf D}_{3t}+\Gamma_1\Delta{\bf Y}_{t-1}+\Pi{\bf Y}_{t-1}+{\bf E}_t.
\end{equation}

\noindent
where $\Pi=\Phi_1+\Phi_2-I_3$, $\Gamma_1=-\Phi_2$, ${\bf c}$ is a vector with constants and ${\bf D}_{it}$ denote the seasonal dummies for trimester $i=1,2,3$. 
The results are displayed in Table \ref{tab:johansen}. 

\begin{table}[!h]
\centering
\caption{FBST and maximum eigenvalue test applied to Finish data of Johansen and Juselius (1990)}\label{tab:johansen}
\begin{tabular}{lccc}
\hline
$H_0$		&\textbf{FBST } 	& \textbf{Max.}	&$p$-value\\ \hline
$r=0$			&0.132  &38.489   &0.0007	\\
$r=1$			&0.994  &26.642   &0.0060	\\
$r=2$			&$\simeq$1  &7.8924   &0.3983	\\ \hline
\end{tabular}
\end{table}

In \cite{johansen1990}, the authors concluded that there is, at least, two cointegration vectors, a conclusion that follows if one adopts a 0.01 significance level, for instance.
Using the asymptotic relationship between $p$-values and e-values for Equation \eqref{eq:finland}, we obtain, for $r=0$, an e-value of 0.998, and, for $r=1$, an e-value of 0.999, corresponding to a 0.01 $p$-value.
These apparently discrepant values for the e-values are due to the high dimensions of the unrestricted ($m=58$) and under $H_0$ ($h=42$ for $r=0$ and $h=43$ for $r=1$) parameter spaces. 
Therefore, under this criterion, the FBST also rejects the null rank and $r=1$ (since 0.132<0.998 and 0.994<0.999, respectively) and does not reject $r=2$, recommending the same action as the maximum eigenvalue test.

\

\noindent{\textbf{Example 4}: \emph{As a final example, we apply the FBST to a US data set discussed in \cite{lucas2000}. The observations have annual periodicity and went from 1900 to 1985. We tested for cointegration between real national income, $M1$ monetary aggregate deflated by the GDP deflator and the commercial papers return rate. The chosen model was a VAR(1) with unrestricted constant. 
The series were used in natural logarithms and the results follow below:}}\\

\begin{table}[!h]
\centering
\caption{FBST and maximum eigenvalue test applied to US data of Lucas (2000)}\label{tab:lucas}
\begin{tabular}{lccc}
\hline
$H_0$		&\textbf{FBST}  	& \textbf{Max}.	&$p$-value\\ \hline
$r=0$			&0.042  &25.334    & 0.0101	\\
$r=1$			&0.996  &4.2507    & 0.8271	\\ \hline
\end{tabular}
\end{table}

Table \ref{tab:lucas} shows that the maximum eigenvalue test rejects $r=0$ and does not reject $r=1$ at a 0.05 significance level.
Once more adopting the asymptotic relationship between $p$-values and e-values for the chosen model, we obtain, for $r=0$, an e-value of 0.247 corresponding to a 0.01 $p$-value.
Thus, under this criterion, the FBST also rejects the null rank and does not reject $r=1$.

\section{Conclusions}

In the past few decades, the econometric literature introduced statistical tests to identify unit roots and cointegration relationships in time series. 
The Bayesian approach applied to these topics advanced considerably after the 1990s, developing interesting alternatives, mostly for unit root testing.  
The (parametric) frequentist tests mentioned here may not be suitable since these procedures rely on the distribution of the test statistic---usually assuming the hypothesis being tested is true---which depend on a particular a statistical model, usually Gaussian.
When the distributions of such statistics cannot be obtained, the procedure is saved by asymptotic results. 
If the researcher considers different statistical models and the available sample is small, the results of the tests may be quite misleading.

The present work reviewed a simple and powerful Bayesian procedure that can be applied to both purposes: unit root and cointegration testing.
We have also shown that the FBST works considerably well even when one uses improper priors, a choice that may preclude the derivation of Bayes Factors, a standard Bayesian procedure in hypotheses testing.   

A long series of articles provided in \cite{stern2020} and the references therein, has showed the versatility and properties of FBST, such as: a. the e-value derivation and computation are straightforward from its general definition; b. it uses absolutely no artificial restrictions like a distinct probability measure on the hypothesis set, induced by some specific parametrization; c. it is in strict compliance with the likelihood principle; d. it can conduct the test with any prior distribution; e. it does not need closed conjectures concerning error distributions, even for small samples; f. it is an exact procedure, since it does rely on asymptotic assumptions; and g. it is invariant with respect to the null hypothesis parametrization and with respect to the parameter space parametrization.
See \cite{stern2013}, p.253 for this property.

To proceed with this research agenda, it would be interesting to perform more simulation studies with the FBST applied to unit root testing for a larger group of parametric and semi-parametric models (likelihoods).  
Another possibility is to include moving average terms in the data generating processes and work with Gaussian and non-Gaussian ARMA models. 
Notice that, given the points made above, these extensions would not impose major problems to the FBST as they would to the frequentist procedures. 
Regarding cointegration, the same extensions may be studied in future works, although the adoption of statistical models outside the Gaussian family would require further efforts to numerically implement the FBST.
We shall also investigate the effect of the prior choice in the estimates of cointegration relations, especially for small samples.

\section*{Acknowledgments}

This work was supported by UFSCar - Federal University of S\~ao Carlos, USP - University of S\~ao Paulo, and UFMS - Federal University of Mato Grosso do Sul. This work was also partially supported by CNPq - the Brazilian National Council of Technological and Scientic Development (grants PQ 307648-2018-4, 302767-2017-7, 301892-2015-6 and 308776-2014-3); and FAPESP - the State of S\~ao Paulo Research Foundation (grants CEPID Shell-RCGI201450279-4; CEPIDCeMEAI 2013-07375-0). The authors are extremely grateful for the support received from their colleagues, collaborators, users and critics in the construction works of their research projects.

\appendix

\section{Computational Resources}\label{app:computer}

The FBST was implemented in all the examples using codes written by the authors in Matlab/Octave programming language.
The results displayed in Tables \ref{tab:ur_results}, \ref{tab:prior}, \ref{tab:during}, \ref{tab:ex_clim}, \ref{tab:johansen}, and \ref{tab:lucas} were obtained using GNU Octave version 4.4.1. The only package required to run the routines was the \texttt{statistics} package (version 1.4.1), necessary to simulate vectors of random variables from the distributions mentioned in the text. 
The codes are briefly described at \url{https://www.ime.usp.br/~jstern/software/}, where they can be freely downloaded. 

The original data sets used in the examples presented in this work can be obtained from the following sources:

\begin{enumerate}

\item Table \ref{tab:ur_results}: fourteen U.S. economic time series used by \cite{schotman1991}. Available at the R library \texttt{urca}, where it is named ``\texttt{npext}''.

\item Example 1: the original series used in \cite{shoeb2009} and \cite{ostergaard2017} are available at \url{https://physionet.org/content/chbmit/1.0.0/}.
The data for the subject analyzed in Example 1 is from file \texttt{chb01\_03.edf}, found inside folder \texttt{chb01}. 
To obtain Tables \ref{tab:prior} and \ref{tab:during}, the data were transformed as described in Example 1.

\item Example 2: the original data set used in \cite{turasie2016} is available at \url{https://climexp.knmi.nl/NCEPData/ghcn_cams_05.nc}, provided by the Global Historical Climatology Network (GHCN)/Climate Anomaly Monitoring System (CAMS).
The data set studied here is the 2 m temperature analysis (0.5$\times$0.5) data, a high resolution (0.5$\times$0.5 degrees in latitude and longitude) global land surface temperature data set covering the period 1949 to near present, in our case May 2020.

\item Example 3: the original data set with four macroeconomic series used by \cite{johansen1990} to estimate the money demand of Finland is available in the R library \texttt{urca} with the name ``\texttt{finland}''.

\item Example 4: the original data used in \cite{lucas2000} can be downloaded from \url{https://www.ime.usp.br/~jstern/software/}.  
\end{enumerate}

\section{Non-Standard Distributions Used in This Article}\label{app:dist}

\subsection{Inverse-Gamma}

The probability density function of the Inverse-Gamma distribution is given 

\[
f_0(x\mid a, b)
= \frac{b^a}{\Gamma(a)}\cdot\left({1\over x}\right)^{a+1}\exp\left(-{b\over x}\right)
\]

\noindent
for $x>0$ and zero, otherwise. 
The parameters $a$ and $b$ are both positive real numbers and $\Gamma$ is the gamma function.

\subsection{Matrix Normal}

The probability density function of the random matrix ${\bf X}$ with dimensions $p\times q$ that follows the matrix normal distribution $MN_{p\times q}(\mathbf{M}, \mathbf{U}, \mathbf{V})$ has the form:

\[
f_1(\mathbf{X}\mid\mathbf{M}, \mathbf{U}, \mathbf{V}) = \frac{\exp\left( -\frac{1}{2} \, \mathrm{tr}\left[ \mathbf{V}^{-1} (\mathbf{X} - \mathbf{M})' \mathbf{U}^{-1} (\mathbf{X} - \mathbf{M}) \right] \right)}{(2\pi)^{pq/2} |\mathbf{V}|^{p/2} |\mathbf{U}|^{q/2}}
\]

\noindent
where $\mathbf{M}\in\mathbb{R}^{p\times q}$, $\mathbf{U}\in\mathbb{R}^{p\times p}$ and $\mathbf{V}\in\mathbb{R}^{q\times q}$, being $\mathbf{U}$ and $\mathbf{V}$ symmetric positive semidefinite matrices.
The matrix normal distribution can be characterized by the multivariate normal distribution as follows: $\mathbf{X} \sim MN_{p\times q}(\mathbf{M}, \mathbf{U}, \mathbf{V})$ if and only if $\mathrm{vec}(\mathbf{X}) \sim N_{pq}(\mathrm{vec}(\mathbf{M}), \mathbf{V} \otimes \mathbf{U})$, where $\otimes$ denotes the Kronecker product and vec the vectorization of $\mathbf{M}$.


\subsection{Inverse-Wishart}

The probability density function of the Inverse-Wishart distribution is

\[
f_2({\mathbf x}\mid\Lambda, \nu) = \frac{\left|{\Lambda}\right|^{\nu/2}}{2^{\nu p/2}\Gamma_p(\frac \nu 2)} \left|\mathbf{x}\right|^{-(\nu+p+1)/2} \exp\left[-\frac{1}{2}\operatorname{tr}(\Lambda\mathbf{x}^{-1})\right]
\]

\noindent
where $\mathbf{x}$ and $\Lambda$ are $p\times p$ positive-definite matrices, and $\Gamma_p$ is the multivariate gamma function.
Notice that we may also write the same density with tr$(\mathbf{x}^{-1}\Lambda)$ inside the exponential function, as would be convenient in our implementation of the Gibbs sampler in Section \ref{sec:impl_coint}.

\section{Heuristic Proof of Johansen's Procedure}\label{app:johansen}
\unskip

The goal of this appendix is to provide a brief heuristic explanation of the procedure, discussed in Section \ref{sec:impl_coint} that finds the maximum of posterior \eqref{eq:pos} subject to the hypothesis that matrix $\Pi$ has reduced rank $r$, $0\leq r\leq n$.
The procedure is based on the algorithm proposed in \cite{johansen1988, johansen1990} to maximize a Gaussian likelihood under the same assumption (reduced rank of matrix $\Pi$).
The formal proof of Johansen's algorithm can be found in \cite{hamilton1994}, chapter 20.
As mentioned in Section \ref{sec:impl_coint}, Johansen's algorithm can be applied to the posterior \eqref{eq:pos} since this distribution is very close to a (multivariate) Gaussian likelihood.

The first step of the algorithm involves ``concentrating'' the posterior, i.e., to assume $\Omega$ and $\Pi$ are given and maximize the posterior with respect to all the other parameters in $\Theta$.
Hence, let $\gamma$ denote the matrix $\eta$ except for matrix $\Pi$, i.e., $\gamma ~ '=\left[ {\bf c} ~ ~ \Phi_0' ~ ~ \Gamma_1' ~ ~ \ldots ~ ~ \Gamma_{p-1}' \right]$.
The concentrated log-posterior, denoted by $\mathcal{M}$, is found by replacing $\gamma$ with $\widehat{\gamma}(\Pi)$ in \eqref{eq:pos}:

\begin{equation}\label{eq:logconc}
\mathcal{M}(\Pi,\Omega\mid{\bf y})=\ln[g(\widehat{\gamma}(\Pi);\Pi,\Omega\mid{\bf y})]=C+ {(T+n+1)\over 2}\ln|\Omega^{-1}|-\left\{-\frac{1}{2}\cdot \textrm{tr}[\Omega^{-1}(\widehat{{\bf U}}-\Pi\widehat{{\bf V}})'(\widehat{{\bf U}}-\Pi\widehat{{\bf V}})]\right\}
\end{equation}

\noindent
where $C$ is a constant that depends on $T$, $n$ and ${\bf y}$.
The strategy behind concentrating the posterior is that, if we can find the values $\widehat{\Omega}$ and $\widehat{\Pi}$ that maximize $\mathcal{M}$, then these same values, along with $\widehat{\gamma}({\widehat{\Pi}})$, will maximize \eqref{eq:pos} under the constraint rank$(\Pi)=r$.
Carrying the concentration on one step further, we can find the value of $\Omega$ that maximizes \eqref{eq:logconc} assuming $\Pi$ known, giving

\[
\widehat{\Omega}(\Pi)={1\over T+n+1}\cdot(\widehat{{\bf U}}-\Pi\widehat{{\bf V}})'(\widehat{{\bf U}}-\Pi\widehat{{\bf V}}).
\]

To evaluate the concentrated log-posterior at $\widehat{\Omega}(\Pi)$, notice that

\[
\text{tr}\left[\widehat{\Omega}(\Pi)^{-1}(\widehat{{\bf U}}-\Pi\widehat{{\bf V}})'(\widehat{{\bf U}}-\Pi\widehat{{\bf V}})\right]=\text{tr}[(T+n+1)I_n]=n(T+n+1)
\]

\noindent
and, therefore, denoting by $\mathcal{M}^*$ this new concentrated log-posterior, we have

\begin{align}
\mathcal{M}^*(\Pi\mid{\bf y})&=C+ {(T+n+1)n\over 2}-{(T+n+1)\over 2}\ln\left|\frac{1}{T+n+1}(\widehat{{\bf U}}-\Pi\widehat{{\bf V}})'(\widehat{{\bf U}}-\Pi\widehat{{\bf V}})\right|\\
&=C+ {(T+n+1)n\over 2}-{(T+n+1)\over 2}\ln\left|\frac{T}{T+n+1}\cdot{1\over T} (\widehat{{\bf U}}-\Pi\widehat{{\bf V}})'(\widehat{{\bf U}}-\Pi\widehat{{\bf V}})\right|\\
&=C+ {(T+n+1)n\over 2}-{(T+n+1)\over 2}\ln\left[\left(\frac{T}{T+n+1}\right)^n\cdot\left|{1\over T} (\widehat{{\bf U}}-\Pi\widehat{{\bf V}})'(\widehat{{\bf U}}-\Pi\widehat{{\bf V}})\right|\right]\\
&=K-{(T+n+1)\over 2}\cdot\ln\left|{1\over T} (\widehat{{\bf U}}-\Pi\widehat{{\bf V}})'(\widehat{{\bf U}}-\Pi\widehat{{\bf V}})\right|\label{eq:conc}
\end{align}

\noindent
where $K$ is a new constant depending only on $T$, $n$ and ${\bf y}$.
Equation \eqref{eq:conc} represents the maximum value one can achieve for the log-posterior for any given matrix $\Pi$.
Thus, maximizing the posterior comes down to choosing $\Pi$ so as to minimize the determinant

\[
\left|{1\over T} (\widehat{{\bf U}}-\Pi\widehat{{\bf V}})'(\widehat{{\bf U}}-\Pi\widehat{{\bf V}})\right|
\]

\noindent
subject to the constraint rank$(\Pi)=r$.
The solution of this problem demands the analysis of the sample covariance matrices of the OLS residuals $\widehat{{\bf U}}$ and $\widehat{{\bf V}}$ and here we only present the final expression for the maximum value achieved for the log-posterior, denoted $\ell^*$ in Section \ref{sec:impl_coint}:

\begin{equation}\label{eq:max_r}
\ell^*=K-{(T+n+1)\over 2}\cdot\ln|\widehat{\Sigma}_{{\bf U}{\bf U}}|-{T+n+1\over 2}\cdot\sum_{i=1}^r \ln(1-\widehat{\lambda}_i).
\end{equation}

Chapter 20 of \cite{hamilton1994} provides the formal derivation of \eqref{eq:max_r}.




\end{document}